\input amstex
\input amsppt.sty
\pagewidth {125mm}
\pageheight {185mm}
\parindent=8mm
\frenchspacing 
\input xy
\xyoption{all}


\NoBlackBoxes
\newdir{ >}{{}*!/-7pt/@{>}}
\define\fib{\ar @{->>} [r]} 
\define\cof{\ar @{ >->}[r]} 
\define\we{\ar [r]^{\sim}}

\define\fls{\Cal L X}
\define\ev {\text{ev}}
\define\si {s^{-1}}
\define\Om {\Omega}
\define\sn #1{(-1)^{#1}}
\define \vp{\varphi}
\define\del{\partial}
\define \cyl{X\sp I}
\define \tot{\widetilde{\otimes}}
\define \ba{\bar A}
\define \bd{\bar d}
\define\tom{\widetilde \Omega}
\define\ssimp{\Cal S_{\bullet}}
\define\im{ \text{Im}\, }
\define\coker{ \text{coker}\, }
\define\bc{\Cal B}

\define\tor{\operatorname{Tor}} 
\define\tD{\widetilde D}
\define\ve{\varepsilon}
\define\bbd{\check D}
\define\wU{\widetilde \Upsilon}
\define\Ups{\Upsilon}
 \define\pth{$p^{\text{th}}$}
 \define\trf{tr\negthinspace f}
\define\Diff{\operatorname{Diff}}
\define\Homeo{\operatorname{Homeo}}
\redefine\H {\operatorname H}
 

\topmatter
\title An Algebraic Model for Mod 2 Topological Cyclic Homology\endtitle
\author Kathryn Hess\endauthor
\date 07.12.04\enddate
\affil Ecole Polytechnique F\'ed\'erale de Lausanne \endaffil
\address Ecole Polytechnique F\'ed\'erale 
de Lausanne (EPFL), Institut de g\'eom\'etrie, alg\`ebre et topologie,  CH-1015 Lausanne, Switzerland\endaddress
\email kathryn.hess\@epfl.ch\endemail

\leftheadtext {K. Hess}
\rightheadtext{Mod 2 Topological Cyclic Homology}

 
 \abstract For any space $X$ with the homotopy type of  simply-connected, finite-type CW-complex, we construct an associative cochain algebra $fls^*(X)$ such that $\H^*(fls^*(X))\cong \H^*(\Cal LX)$ as algebras, where $\Cal LX$ denotes the free loop space on $X$.  Under additional conditions on $X$, e.g., when $X$ is a wedge of spheres, we define a cochain complex $hos^*(X)$ by twisting together $fls^*(X)$ and $\H^*(BS^1)$ and prove that $\H^*(hos^*(X))\cong \H^*(\Cal LX_{hS^1})$ as graded modules.  We also show that there is a cochain  map from $fls^*(X)$ to itself that is a good model of the \pth -power operation on $\Cal LX$.  Finally, we define $tc^*(X)$ to be the mapping cone of the composite of the projection map from $hos^*(X)$ to $fls^*(X)$ with the model of the \pth-power map (for $p=2$), so that the mod 2 spectrum cohomology  of $TC(X;2)$ is isomorphic to $\H^*(tc^*(X)\otimes \Bbb F_2)$.  We conclude by calculating $\H^*(TC(S^{2n+1};2);\Bbb F_2)$.
 \endabstract \endtopmatter
\document


\head Preface\endhead

\noindent B\"okstedt, Hsiang and Madsen introduced {\sl topological cyclic homology (TC)} as a
topological version of Connes' cyclic homology in \cite {BHM}.  The topological cyclic
homology of a space $X$ at a prime $p$, denoted $TC(X;p)$, is a spectrum that is the target of the 
{\sl cyclotomic trace map}
$$\xymatrix@1{Trc :A(X)\ar [r]&TC(X;p),}$$
the source of which is Waldhausen's algebraic $K$-theory spectrum of $X$, to which $TC(X;p)$ provides a useful approximation. The cyclotomic trace map is analogous to the Dennis trace map $K_*(A)@>>> HH_*(A)$.  Very little is known about the $TC(X;p)$ when $X$ is not a singleton.  

Waldhausen's algebraic $K$-theory itself approximates the smooth and topological Whitehead spectra $Wh^d(X)$ and $Wh^t(X)$. There are natural cofiber sequences of spectra
$$
\xymatrix{ \Sigma ^\infty X_+ \ar [r]^{\eta_X}& A(X) \ar [r]& Wh^d(X)}
$$
and
$$
\xymatrix{A(*) \wedge X_+ \ar [r]^(0.55){a^A_X}& A(X) \ar [r]& Wh^t(X)}
$$
respectively, where $\eta_X$ is the unit map and
$a^A_X$ is the $A$-theory assembly map \cite{W1, 3.3.1}. Here, $Y_+$ denotes the  space $Y$ with an extra basepoint,  and $\Sigma ^\infty Y$ denotes the suspension spectrum of $Y$. 

By Waldhausen's stable parametrized $h$-cobordism theorem \cite{W2},
there are homotopy equivalences
$$
\Cal H^d(M) \simeq \Omega \Omega^\infty Wh^d(M)
$$
and
$$
\Cal H^t(M) \simeq \Omega \Omega^\infty Wh^t(M)
$$
when $M$ is a smooth, respectively topological, compact manifold.  Here $\Cal
H^d(M)$ is the stable smooth $h$-cobordism space of $M$, which in a stable
range carries information about the homotopy type of the topological
group $\Diff(M)$ of self-diffeomorphisms of $M$.   Likewise $\Cal H^t(M)$
is the stable topological $h$-cobordism space of $M$, which in a stable
range carries information about the topological group $\Homeo(M)$ of
self-homeomorphisms of $M$.

Any information we may obtain about the topological cyclic homology of a smooth or topological manifold will therefore give us some indication of the nature of its stable $h$-cobordism space.

For further explanation of the role of trace maps in $K$-theory and (topological) cyclic and Hochschild homology, as well as about Whitehead spectra, we refer the reader to \cite {Be} and \cite {R}.  

The goal of this chapter, as well as of the minicourse upon which the chapter is based, is to construct  a cochain complex $tc^*(X)$ such that $\H^*(tc^*(X)\otimes \Bbb F_p)$ is isomorphic to the mod $p$ spectrum cohomology of $TC(X;p)$.  For reasons of ease of notation and computation, we will limit ourselves to $p=2$ in this article. 

There are several equivalent definitions of $TC(X;p)$.  The definition that is best suited to algebraic modeling can be stated as follows \cite {BHM}.  Let $\Cal LX$ be the free loop space on $X$, i.e., the space of unbased maps from the circle $S^1$ into $X$, which admits a natural $S^1$-action, by rotation of loops.  Let $\Cal LX _{hS^1}=ES^1\underset S^1\to\times \Cal L X$ denote the homotopy orbit space of this action.  Let $\xymatrix@1{\lambda ^p:\Cal LX\ar [r]&\Cal LX}$ denote the \pth-power map, defined by $\lambda ^p(\ell)(z)=\ell (z^p)$ for all $\ell \in \Cal L X$ and all $z\in S^1$.   There is a homotopy pullback of spectra
$$\xymatrix{TC(X;p)\ar [r]\ar [d]&\Sigma ^\infty \Cal LX_+\ar [d]^{\Sigma ^\infty (Id-\lambda ^p)}\\
\Sigma ^\infty\bigl(\Sigma (\Cal LX_{hS^1})_+\bigr)\ar [r]^(0.6){\trf_{S^1}}&\Sigma ^\infty\Cal L X_+}$$
where $\trf_{S^1}$ is the $S^1$-transfer map associated to the principal $S^1$-bundle
$$\xymatrix{ES^1\times \Cal L X\ar [r]&\Cal LX_{hS^1}.}$$
It is therefore clear that $TC(X;p)$ is the homotopy fiber of the composition
$$\xymatrix{\Sigma ^\infty \Cal LX_+\ar [rr]^{\Sigma ^\infty (Id-\lambda ^p)}&&\Sigma ^\infty \Cal LX_+\ar [r]^(0.4)\iota& \operatorname{hocofib} (\trf_{S^1})}.$$

Motivated by this characterization of $TC(X;p)$, we apply the following method to constructing $tc^*(X)$.  We first define an associative cochain algebra $fls^*(X)$ together with a cochain map 
$$\xymatrix{\Ups:fls ^*(X)\ar [r]&CU^*\Cal LX}$$
inducing an isomorphism of algebras in cohomology, where $CU^*$ denotes the reduced cubical cochains.  We then twist together $fls^*(X)$ and $\H^*(BS^1)$, obtaining a new cochain complex $hos^*(X)$ that fits into a commuting diagram
$$\xymatrix{hos^*(X)\ar [d]^{\overline\Ups}\ar [r]^\pi &fls^*(X)\ar [d]^\Ups\\
		   CU^*(\Cal L X_{hS^1})\ar [r]^{CU^*c}&CU^*(\Cal LX)}$$
where $\pi$ is the projection map, $\xymatrix@1{c:\Cal LX\ar [r]&\Cal L X_{hS^1}}$ is the map induced by the inclusion $\xymatrix@1{\Cal LX\ar [r]&ES^1\times \Cal LX}$ and $\overline\Ups$ induces an isomorphism in cohomology.  The projection map $\xymatrix@1{\pi:hos^*(X)\ar [r]&fls^*(X)}$ is then a model for the inclusion $$\xymatrix@1{\Sigma ^\infty \Cal LX_+\ar [r]^(0.4)\iota& \operatorname{hocofib} (\trf_{S^1})}.$$  Finally, we define a cochain map $\xymatrix@1{\frak l^p:fls^*(X)\ar [r]&fls^*(X)}$ such that $\Ups\circ \frak l^p$ and $CU^*\lambda ^p\circ \Ups$ are chain homotopic. Thus $Id-\frak l^p$ is a model for $\Sigma ^\infty (Id-\lambda ^p)$.  As explained carefully and  in detail in \cite {HR}, if we set $tc^*(X)$ equal to the mapping cone of the composition $(Id-\frak l^p)\pi$, then
$\H^*(tc^*(X)\otimes \Bbb F_p)\cong \H^*(TC(X;p);\Bbb F_p)$, as desired.  Here we remark only that there is a Thom isomorphism involved in the
identification of the cohomology of the homotopy cofiber of the
$S^1$-transfer map and the cohomology of the $S^1$-homotopy orbits of the free
loop space.  Furthermore, the fact that the projection map $\pi$ is a model for
the inclusion $\iota$ requires an analysis of these transfers and
isomorphisms

The article is organized as follows. We begin in section 0 by reminding the reader of certain algebraic and topological notions and constructions. In section 1 we study free loop spaces and their algebraic models, beginning by defining a simplicial set that models $\Cal L X$, which we then apply to constructing $fls^*(X)$ via a refined version of methods from \cite {DH1-4}.  Section 2 is devoted to the study of homotopy orbit spaces of circle actions.  We first treat the general case, twisting together $\H^*(BS^1)$ and $CU^*Y$ to obtain a large but attractive cochain complex $HOS^*(Y)$ for calculating $\H^*(Y_{hS^1})$, when $Y$ is any $S^1$-space.  Specializing to the case $Y=\Cal LX$, we show how to twist together $\H^*(BS^1)$ and $fls^*(X)$ to build $hos^*(X)$,  so that we obtain a complex equivalent to $HOS^*(\Cal LX)$.  Finally, in section 3 we define and study our model for the \pth-power map (for $p=2$) , then apply it, together with the results of the preceding chapters, to the construction of $tc^*(X)$.  We conclude by applying our model to the calculation of the mod 2 spectrum cohomology of $TC(S^{2n+1};2)$.

\remark{Remark} In these notes, complete proofs are provided only of those results that have not yet appeared elsewhere and that are due to the author. Furthermore, some results that have yet to be published are not proved completely here, if the complete proof is excessively technical. We hope in such cases to have provided enough detail to convince the reader of the truth of the statement.   The reader who is curious about the details is refered to articles that should appear soon.\endremark

The author would like to thank David Chataur, John Rognes and J\'er\^ome Scherer for their helpful comments on earlier versions of this chapter.  Warm thanks are also due to David Chataur, Jos\'e-Luis Rodrigues and J\'er\^ome Scherer for their splendid organization of the Almer\'\i a summer school on string topology.


\head 0. Preliminaries\endhead

\noindent We begin here by recalling certain elementary definitions and constructions and fixing our basic notation and terminology. We then remind the reader of the construction of the canonical, enriched Adams-Hilton model of a simplicial set, which is the input data for our free loop space model.  We conclude this section with a description of our general method for constructing algebraic models of fiber squares, which we then apply in section 1 to building our free loop space model.

\subhead 0.1 Elementary definitions, terminology and notation\endsubhead

\noindent Throughout this paper we work over $\Bbb Z$, the ring of 
integers, unless stated otherwise.

Given chain complexes $(V,d)$ and $(W,d)$, the notation
$f:(V,d)@>\simeq >>(W,d)$ indicates that $f$ induces an isomorphism in homology. 
In this case we refer to $f$ as a {\sl quasi-isomorphism}.

If $V=\bigoplus \sb {i\in \Bbb Z} V\sb i$ is a graded module, then
 $\si V$ and $sV$ denote the graded modules with, respectively, 
$(\si V)\sb i \cong V\sb {i+1}$ and $(sV)\sb i \cong V\sb {i-1}$.  Given a homogeneous element $v$ in
$V$, we write $\si v$ and $sv$ for the corresponding elements of $\si V$ and $sV$.  If the gradings are written as upper indices, i.e., $V=\bigoplus \sb {i\in \Bbb Z} V^ i$, then $(\si V)^ i \cong V^{i-1}$ and $(sV)^ i \cong V^ {i+1}$.

Dualization is indicated throughout the paper by a $\sharp$ as superscript.  The degree of an element $x$ in a graded module is denoted $|x|$, unless it is used as an exponent, in which case the bars may be dropped.

A graded $\Bbb Z$-module $V=\bigoplus \sb {i\in \Bbb Z} V\sb i$ is {\sl connected} 
if $V_{<0}=0$ and $V_{0}\cong \Bbb Z$.  It is {\sl simply connected} if, in 
addition, $V_{1}=0$.  We write $V_{+}$ for $V_{>0}$.
Let $V$ be a positively-graded, free $\Bbb Z$-module.  The free associative algebra on $V$ is denoted 
$TV$, i.e., 
$$TV\cong\Bbb Z\oplus V\oplus (V\otimes V)\oplus (V\otimes V\otimes V)\oplus\cdots .$$
A typical basis element of $TV$ is denoted $v\sb 1\cdots v\sb n$, i.e., we drop the
tensors from the notation.  The product on $TV$ is then defined by
$$\mu (u\sb 1\cdots u\sb m\otimes v\sb 1\cdots v\sb n)=u\sb 1\cdots u\sb m v\sb
1\cdots v\sb n.$$

The cofree, coassociative coalgebra on $V$, denoted
$\bot V$ in this article, is isomorphic as a graded $\Bbb Z$-module
to $TV$.  We write $\bot \sp n V=\bigotimes \sp n V$,
of which a typical basis element is denoted $v\sb 1|\cdots |v\sb
n$. 
The coproduct on $\bot V$ is then defined in the obvious manner
by
$$\split
\Delta (v\sb 1|\cdots| v\sb n)&=v\sb 1|\cdots| v\sb n\otimes 1 +1
\otimes v\sb 1|\cdots |v\sb n\\
&\quad +\sum \sb {i=1}\sp {n-1}v\sb 1|\cdots |v\sb
i\otimes v\sb {i+1}|\cdots |v\sb n.
\endsplit $$

Let $(C,d)$ be a simply-connected (co)chain 
coalgebra with reduced coproduct $\overline \Delta$.  The {\sl cobar construction} on $(C,d)$, denoted $\Om 
(C,d)$, is the (co)chain algebra $(T\si (C_{+}),d_{\Om})$, where 
$d_{\Om}=-\si ds+(\si \otimes \si)\overline \Delta s$ on generators. 

 Let $(A,d)$ be a connected chain algebra or a
simply-connected cochain algebra over
$R$, and let $\bar A$ be the component of $A$ of positive degree. The
{\sl bar construction} on
$(A,d)$, denoted
$\Cal B(A,d)$, is a differential graded coalgebra $(\bot (s\bar A),
D\sb {\Cal B})$.  Let $(D\sb
{\Cal B})\sb 1$ denote the linear part of the differential, i.e.,
$(D\sb {\Cal B})\sb 1=\pi D\sb {\Cal B}$, where $\pi :\bot V@>>>
V$ is the natural projection.  The linear part of $D\sb {\Cal B}$
specifies the entire differential and is given by 
$$(D\sb {\Cal
B})\sb 1(sa\sb 1|\cdots |sa\sb n)=\left\{\aligned  -s(da\sb 1)\quad
&\text{if}\quad n=1\\ \sn {a\sb 1+1}s(a\sb 1\cdot a\sb 2)\quad
&\text{if}\quad n=2\\ 0\quad &\text{if}\quad n>2.
\endaligned\right.$$

\definition {Definition}  Let $f,g:(A,d)@>>>(B,d)$ be two maps of 
chain (respectively, cochain) algebras.  An {\sl $(f,g)$-derivation homotopy} is a linear map $\vp 
:A@>>>B$ of degree $+1$ (respectively, $-1$) such that $d\vp+\vp d= f-g$ and $\vp \mu 
=\mu (\vp \otimes g+f \otimes \vp)$, where $\mu$ denotes the 
multiplication on $A$ and $B$. \enddefinition

If $f$ and $g$ are maps of (co)chain coalgebras, there is an obvious 
dual definition of an {\sl $(f,g)$-coderivation homotopy}.

We often apply Einstein's summation convention in this chapter. When an index appears as both a subscript and a superscript in an expression, it is understood that we sum over that index.  For example, given an element $c$ of a coalgebra $(C, \Delta)$, the notation $\Delta (c) = c_i\otimes c^i$ means $\Delta (c) =\sum _{i\in I} c_i\otimes c^i$.

Another convention used consistently throughout this chapter is the Koszul sign convention for commuting elements  of a graded module or for commuting a morphism of graded modules past an element of the source module.  For example,  if $V$ and $W$ are graded algebras and $v\otimes w, v'\otimes w'\in V\otimes W$, then $$(v\otimes w)\cdot (v'\otimes w')=(-1)^{|w|\cdot |v'|}vv'\otimes ww'.$$ Futhermore, if $f:V@>>>V'$ and $g:W@>>>W'$ are morphisms of graded modules, then for all $v\otimes w\in V\otimes W$, 
$$(f\otimes g)(v\otimes w)=(-1)^{|g|\cdot |v|} f(v)\otimes g(w).$$
The source of the Koszul sign convention is the definition of the twisting isomorphism
$$\tau: V\otimes W@>>>W\otimes V: v\otimes w\mapsto (-1)^{|v|\cdot |w|}w\otimes v.$$

We assume throughout this chapter that the reader is familiar with the elements of the theory of simplicial sets and of model categories.    We recall here only a few very basic definitions, essentially to fix notation and terminology, and refer the reader to, e.g., \cite {May} and \cite {GJ} for simplicial theory and to \cite {Ho}, \cite {DS} and \cite {H} for model category theory.

\definition {Definition} Let $K$ be a simplicial set, and let $\Cal 
F_{ab}$ denote the free abelian group functor.  For all $n>0$, let 
$DK_{n}=\cup _{i=0}^{n-1} s_{i}(K_{n-1})$, the set of degenerate 
$n$-simplices of $K$. The {\sl 
normalized chain complex} on $K$, denoted $C_*(K)$, is given by
$$C_{n}(K)=\Cal F_{ab} (K_{n})/\Cal F_{ab}(DK_{n}).$$
Given a map of simplicial sets $f:K@>>>L$, the induced map of 
normalized chain complexes is denoted $C_*f$.
\enddefinition

Recall that $H_*(C_*(K))\cong H_*(|K|)$ as graded coalgebras, where $|K|$ denotes the geometric realization of $K$.

\definition {Definition} Let $K$ be a reduced simplicial set, and let $\Cal F$ denote the 
free group functor. The {\sl loop group} $GK$ on $K$ is the simplicial group such that $(GK)_{n}=\Cal F 
(K_{n+1}\smallsetminus Im s_{0})$, with faces and degeneracies 
specified by
$$\split 
\del _{0}\bar x&=(\;\overline {\del _{0}x}\;)^{-1}\overline {\del 
_{1}x}\\
\del _{i}\bar x&=\overline {\del 
_{i+1}x}\quad\text {for all $i>0$}\\
s_{i}\bar x&=\overline {s_{i+1}x}\quad \text {for all $i\geq 0$}
\endsplit$$
where $\bar x$ denotes the class in $(GK)_{n}$ of $x\in K_{n+1}$.
\enddefinition

Recall that $H_*(GK)\cong H_*(\Om |K|)$ as graded Hopf algebras.

In any model category we use the notation $\xymatrix@1{\cof&}$ for cofibrations, $\xymatrix@1{\fib&}$ for fibrations and $\xymatrix@1{\we&}$ for weak equivalences.

\subhead 0.2 The canonical, enriched Adams-Hilton model\endsubhead

\noindent We recall in this section the construction given in \cite {HPST} of the 
canonical, enriched Adams-Hilton model of a $1$-reduced simplicial set 
$K$, upon which our free loop space model construction is based.  We begin by reminding the reader of 
the theories that are essential to this construction.  We first sketch briefly the 
classical and crucial theory of twisting cochains, which goes back to work of E. Brown \cite {Br}.  We then outline 
the theory of strongly homotopy coalgebra maps. We conclude this section  by presenting the canonical 
Adams-Hilton model. 

\subsubhead Twisting cocochains\endsubsubhead \definition{Definition} Let $(C,d)$ be a chain coalgebra with 
coproduct $\Delta$, and let 
$(A,d)$ be a chain algebra with product $\mu$.  A {\sl twisting cochain} from $(C,d)$ 
to $(A,d)$ is a degree $-1$ map $t:C@>>>A$ of graded modules such 
that
$$dt+td=\mu (t\otimes t)\Delta.$$
\enddefinition

The definition of a twisting 
cochain $t:C@>>>A$ is formulated precisely so that the following two 
constructions work smoothly. First, let $(A,d)\otimes 
_{t}(C,d)=(A\otimes C, D_{t}),$ where $D_{t}=d\otimes 1_{C}¥+1_{A}¥\otimes d- 
(\mu\otimes 1_{C}¥)(1_{A}¥\otimes t\otimes 1_{C}¥)(1_{A}¥\otimes \Delta)$.  It is easy to 
see that $D_{t}^2=0$, so that $(A,d)\otimes 
_{t}(C,d)$ is a chain complex, which extends $(A,d)$, i.e., of which $(A,d)$ is subcomplex. 
Second, if $C$ 
is connected, let 
$\tilde t:T\si C_{+}¥@>>>A$ be the algebra map given by $\tilde t(\si 
c)=t(c)$.  Then $\tilde t$ is in fact a chain algebra map $\tilde 
t:\Om (C,d)@>>>(A,d)$.  It is equally clear that any algebra map 
$\theta :\Om (C,d)@>>>(A,d)$ gives rise to a twisting cochain via the 
composition
$$C_{+}¥@>\si>>\si C_{+}\hookrightarrow T\si C_{+}@>\theta >> A.$$
Furthermore, the complex $(A,d)\otimes 
_{t}(C,d)$ is acyclic if and only if $\tilde t$ is a 
quasi-isomorphism.

The twisting cochain associated to the cobar construction is a 
fundamental example of this notion.  Let $(C,d,\Delta )$ be a 
simply-connected chain coalgebra.  Consider the linear map
$$t_{\Om}:C@>>>\Om C:c@>>>\si c.$$
It is an easy exercise to show that $t_{\Om}$ is a twisting cochain and 
that $\tilde t_{\Om}=1_{\Om C}$.  Thus, in particular, $(\Om C,d)\otimes 
_{t_{\Om}¥}(C,d)$ is acyclic; this is the well-known acyclic cobar 
construction.
 
\subsubhead Strongly homotopy coalgebra and comodule maps\endsubsubhead
\noindent In \cite {GM} Gugenheim and Munkholm showed that $\operatorname {Cotor}$ was 
natural with respect to a wider class of morphisms than the 
usual morphisms of chain  coalgebras. Given two 
chain coalgebras $(C,d,\Delta)$ and $(C', d', 
\Delta ')$, a {\sl strongly homotopy coalgebra (SHC) map} $f: (C,d,\Delta )\Rightarrow (C',d', \Delta ')$ is a chain map 
$f:(C,d)@>>>(C',d')$ together with a family of $\Bbb Z$-linear maps
$$\frak F (f)=\{ F_{k}:C@>>> (C')^{\otimes k}\mid \deg F_{k} = k-1,\; k\geq 1\}$$
satisfying
\roster 
\item $F_{1}=f$ and
\item  for all $k\geq 2$ 
$$\multline
F_{k}d-\sum _{i+j=k-1}(-1)^{j}(1_{C'}^{\otimes i}\otimes d'\otimes 1_{C'}^{\otimes 
j})F_{k} \\=\sum _{i+j=k}(-1)^{j}(F_{i}\otimes F_{j})\Delta +\sum 
_{i+j=k-2}(-1)^{j}(1_{C'}^{\otimes i}\otimes \Delta '\otimes 1_{C'}^{\otimes 
j})F_{k-1}.\endmultline$$
\endroster

We call $\frak F(f)$ an {\sl SHC family} for $f$.

An SHC map is thus a coalgebra map, up to an infinite family of 
homotopies.  In particular if $f$ is a map of chain coalgebras, then 
it can be seen as an SHC map, with $F_{k}=0$ for all $k>1$.  Furthermore if  $f: 
(C,d,\Delta )\Rightarrow (C',d', \Delta ')$ is an SHC map and $g: 
(C',d',\Delta ')@>>> (C'',d'', \Delta '')$ is a strict coalgebra map, 
then $gf$ is an SHC map, where $\frak F(gf)=\{ g^{\otimes k}F_{k}\mid 
k\geq 1\}$. 

Observe that the existence of $\frak F(f)$ is equivalent to the 
existence of a chain algebra map $\tom f:\Om (C,d)@>>>\Om (C',d')$ 
such that $\tom f (\si c) - \si f(c)\in T^{\geq 2}\si C'_{+}$.  
Given $\frak F(f)$, we can define $\tom f$ by setting
$$\tom f(\si c)=\sum _{k\geq 1}(\si )^{\otimes k}F_{k}(c)$$
and extending to a map of algebras.
Condition (2) above then implies that $\tom f$ is a differential map 
as well. 

Note that if $f$ is a strict coalgebra map, seen as an SHC map with trivial SHC family, then $\tom f=\Om 
f$.  More generally, if  $f: 
(C,d,\Delta )\Rightarrow (C',d', \Delta ')$ is an SHC map and $g: 
(C',d',\Delta ')@>>> (C'',d'', \Delta '')$ is a strict coalgebra map, seen as an SHC map with trivial SHC family, 
then there is an SHC family for $gf$ such that  $\tom (gf)=\Om g \circ \tom f$.

Similarly, given $\tom f$, we can define $F_{k}$ via the composition
$$C_{+}¥@>\si >>\si C_{+}\hookrightarrow T\si C_{+} @>\tom f>>T\si 
C'_{+}@> proj>>\bigl (\si C'_{+}\bigr ) ^{\otimes k} @>s^{\otimes k}>> (C') ^{\otimes k}.$$

Gugenheim and Munkholm proved in \cite {GM} that the usual simplicial Alexander-Whitney map
$$\xymatrix@1{f_{K,L}:C_*(K\times L)\ar [r]&C_*(K)\otimes C_*(L)}$$
defined by $f_{K,L}(x,y)=\sum _{i=0}^n\del _{i+1}\cdots\del _nx\otimes \del _0^iy$ is naturally an SHC map.

 \subsubhead The canonical Adams-Hilton model\endsubsubhead
\noindent For every pair of simply-connected chain coalgebras $(C,d)$ and $(C',d')$, 
Milgram proved that there is a quasi-isomorphism of chain algebras
$$\rho:\Om \bigl( (C,d)\otimes (C',d')\bigr )@>>>\Om (C,d)\otimes \Om 
(C',d')\tag 0.2.1$$
specified by $\rho\bigl( \si (x\otimes 1)\bigr)=\si x$, $\rho\bigl( \si (1\otimes 
y)\bigr)=\si y$ and $\rho\bigl( \si (x\otimes y)\bigr)=0$ for all $x\in 
C_{+}$ and $y\in C'_{+}$ \cite {Mi}.

In \cite {S} Szczarba gave an explicit formula for a 
natural transformation of functors from simplicial sets to chain algebras
$$\theta: \Om C_*(-)@>>>C_*(G(-))$$
such that $\theta_{K}: \Om C_*(K)@>>>C_*(GK)$ is a quasi-isomorphism of 
chain algebras for every $1$-reduced simplicial set $K$.  Since 
$C_*(GK)$ is in fact a chain Hopf algebra, it is reasonable to ask 
whether $\Om _*C(K)$ can be endowed with a coproduct with respect to 
which $\theta _{K}$ is a quasi-isomorphism of chain Hopf algebras. 

Let $\psi : \Om C_*(-)@>>>\Om C_*(-)\otimes \Om C_*(-)$ denote the 
natural transformation given for each $1$-reduced simplicial set $K$ by the composition
$$\Om C_*(K)@>\Om (\Delta _{K})_{\sharp}>>\Om C_*(K\times K)@>\tom 
f_{K,K}>>\Om \bigl (C_*(K)\otimes C_*(K)\bigr )@>\rho>>\Om C_*(K)\otimes \Om 
C_*(K).$$
The coproduct $\psi _{K}:\Om C_*(K)@>>>\Om C_*(K)\otimes \Om C_*(K)$ is 
called the {\sl Alexander-Whitney (A-W) cobar diagonal}.  In \cite {HPST} Hess, Parent, Scott and Tonks proved that for all $1$-reduced 
$K$, the Alexander-Whitney cobar diagonal is strictly coassociative 
and cocommutative up to derivation homotopy, which we call $\Theta$.  They established 
furthermore that Szczarba's equivalence $\theta_{K}$ is an SHC map with respect 
to $\psi _{K}$ and the usual coproduct on $C_*(GK)$.  

In \cite {B} Baues provided a purely combinatorial definition of strictly coassociative coproduct and of a derivation homotopy for cocommutativity on 
$\Om C_*(K)$ for any $1$-reduced simplicial set $K$, but without giving a map from $\Om C_*(K)$ to $C_*(GK)$.  In \cite {HPST} it is shown that the Alexander-Whitney cobar diagonal is the same as Baues's coproduct, which implies that 
$$\im \overline \psi _K\subseteq T^{\geq 1}\si C_+K\otimes \si C_+K,$$
where $\overline\psi_K$ is the reduced coproduct.

Henceforth we refer to $\theta_{K}:\Om C_*(K)@>>>C_*(GK)$ as the {\sl 
canonical Adams-Hilton model} and to $\psi _{K}:\Om C_*(K)@>>>\Om C_*(K)\otimes 
\Om C_*(K)$ as its {\sl canonical enrichment}.

\subhead 0.3 Noncommutative algebraic models of fiber squares\endsubhead

\noindent We review in this section the bare essentials of noncommutative modeling of fiber squares, as developed in \cite {DH1}.  Note that this theoretical framework is highly
analogous to the theory of KS-extensions in rational homotopy theory.  See \cite {DH1}and
\cite {FHT} for more details.     

We first define the classes of morphisms with which we work throughout the remainder
of this article.

\definition{Definition} Let $(B,d)$ and $(C,d)$ be bimodules over an associative cochain algebra $(A,d)$. A cochain map $\xymatrix@1{f:(B,d)\ar[r]& (C,d)}$  is a {\sl quasi-bimodule map}
if $H\sp * f$ is a map of $\H^*(A,d)$-bimodules. If  $(A,d)$ and $(A',d')$ are associative cochain algebras, then a cochain map $\xymatrix@1{f:(A,d)\ar[r]& (A',d')}$  is a {\sl quasi-algebra map}
if $H\sp * f$ is a map of algebras.
\enddefinition   

Noncommutative cochain algebra 
models of topological spaces are defined in terms of quasi-algebra maps.

\definition{Definition}Let $X$ be a topological space.  An (integral) {\sl noncommutative model} of
$X$ consists of an associative cochain algebra over $\Bbb Z$, $(A,d)$, together with a
quasi-algebra quasi-isomorphism $$\xymatrix@1{\alpha :(A,d)\ar [r]^{\simeq} &C\sp *(X)},$$    
where $\alpha $ is called a {\sl model morphism}.
\enddefinition

Of course, we must also define what it means to model a continuous map, if we wish to model
pull-backs of fibrations. 

\definition{Definition}Let $f:Y@>>>X$ be a continuous map. A {\sl noncommutative 
model} of $f$ consists of a commuting diagram
$$\xymatrix
{(A,d) \ar [r]^\vp \ar [d]^\alpha_\simeq&(B,d)\ar [d]^\beta_\simeq\\
 C\sp *X \ar [r]^{C\sp *f}&   C\sp *Y }$$
in which $\alpha $ and $\beta $ are model morphisms, and $\vp $ is a
quasi-algebra map.
\enddefinition

\remark{Remarks}\roster
\item In most applications of noncommutative models, $\alpha $ and
$\vp$ are strict morphisms of algebras, while it is usually impossible for $\beta $ to be a strict morphism of
algebras. 

\item It is often difficult to use a noncommutative model for
constructions or calculations, unless $A$ is a free algebra.   It is easy to see, however,
that every space $X$ possesses such a noncommutative model.
\endroster\endremark\medskip

There is a special class of strict algebra maps, known as {\sl twisted algebra
extensions}, that are used for modeling topological fibrations.  Roughly speaking, a
twisted algebra extension of one algebra by another is a tensor product of the two
algebras in which both the differential and multiplication are perturbed from the usual
tensor-product differential and multiplication.

\definition{Definition}Let $(A,d)$ and $(B,d)$ be a cochain algebra and a cochain complex over $\Bbb Z$, respectively.  A {\sl twisted
bimodule extension}  of
$(A,d)$ by $(B,d)$ is an $(A,d)$-bimodule $(C,D)$ such that 
\roster
\item $C\cong A\otimes B$ as graded modules;
\item the right action of $A$ on $C$ is free, i.e., $(a\otimes b)\cdot a'=
(-1)\sp {a'b}aa'\otimes b$; 
\item the left action of $A$ on $C$ commutes
with the right action, i.e.,
$$(a\cdot c)\cdot a'=a\cdot (c\cdot a')$$
for all $a,a'\in A$ and $c\in C$, and satisfies 
$$a\cdot(1\otimes b)-a\otimes b\in A\sp +\otimes B\sp {<b}$$ for all $a$ in
$A$ and
$b$ in $B$; and 
\item the inclusion map $(A,d)@>>>(C,D)$ and the projection map $(C,D)@>>>(B,d)$ are
both $(A,d)$-bimodule maps, where $(B,d)$ is considered with the trivial $(A,d)$-bimodule structure, so that, in particular,
$$D(1\otimes b)-1\otimes db\in A\sp +\otimes B$$
for all $b,b'$ in $B$.
\endroster
If $(B,d)$ is a cochain algebra, then a twisted bimodule extension $(C,D)$ of $(A,d)$ by $(B,d)$ is a {\sl twisted algebra extension} if the bimodule structure of $(C,D)$ extends to a full algebra structure such that the inclusion and projection maps above are maps of cochain algebras.  In particular,
$$(1\otimes b)(1\otimes b')-1\otimes bb'\in A\sp +\otimes B$$
for all $b,b'$ in $B$.
\enddefinition

\remark{Notation}We write $(A,d)\tot (B,d)$ to denote a
twisted bimodule extension of $(A,d)$ by $(B,d)$ and $(A,d)\odot (B,d)$ to denote a twisted algebra extension.
\endremark
\medskip

The proposition below, which is the noncommutative analogue of
a well-known result concerning KS-extensions, states that twisted algebra extensions
have the left lifting property with respect to surjective quasi-algebra
quasi-isomorphisms.  Since it is natural to think of surjective
quasi-algebra morphisms as fibrations of cochain algebras, Proposition 0.3.1
implies that we can think of twisted extensions as cofibrations. 
In other words, twisted algebra extensions are plausible models of
topological fibrations, since the cochain functor is contravariant.    
 
\proclaim{Proposition 0.3.1}Let $\iota :(A,d)@>>>(A,d)\odot(B,d)$ be a twisted
algebra extension. Given a commuting diagram 
$$\xymatrix{
(A,d)\ar[r]^f \ar [d]^\iota&( C,d)\ar@{->>}[d]^p_\simeq\\
(A,d)\odot (B,d)\ar[r]^(0.6)g&(E, d)}$$
in which $f$ is a right $(A,d)$-module map, $p$ is a surjective quasi-algebra
quasi-isomorphism, and
$g$ is a quasi-algebra map, there exists a quasi-algebra
map,
$$\xymatrix @1{h: (A,d)\odot (B,d)\ar [r]&(C,d)}$$
which is a right $(A,d)$-module map, as well as a lift of $g$ through $p$ and an extension of $f$, i.e., $ph=g$ and $h\iota=f$ .  
\endproclaim

This proposition is a simplified version of a result that first appeared in \cite {DH1}, but that we do not need in its full generality here. 

\demo{Proof} Since $(A,d)\odot (B,d)$ is semifree as a right $(A,d)$-module, the lift $h$ exists as a map of right $(A,d)$-modules.  In cohomology $\H^*p\H^*h=\H^*g$, which implies that $\H^*h=(\H^*p)^{-1}\H^*g$, since $\H^*p$ is an isomorphism.  Hence $\H^*h$ is an algebra map, as it is a composition of algebra maps.\qed\enddemo

Let us see how to model pull-backs of fibrations in this context.   Consider a pull-back
square of simply-connected spaces  
$$\xymatrix{
E\underset B\to{\times}X \ar [r]^{\bar f}\ar [d]^{\bar q}&E\ar [d]^q\\
X \ar [r]^f&B}$$
in which $q$ is a fibration and $f$ an arbitrary continuous map. Suppose that 
$$\xymatrix{(A,d)\ar[r]^{\vp}\ar [d]^\alpha_\simeq&(\ba ,\bd)\ar [d]^\gamma_\simeq\\
C\sp *B\ar [r]^{C\sp *f}&C\sp *X}\tag 0.3.1$$
and
$$\xymatrix{
(A,d)\ar [r]^(0.4)\iota\ar [d] ^\alpha_\simeq&(A,d)\odot (C,e)\ar [d]^\beta _\simeq\\
C\sp *B\ar [r]^{C\sp *q}&C\sp *E}\tag 0.3.2$$
are noncommutative models of $f$ and $q$, where $\iota $ is
a twisted algebra extension of $(A,d)$.  We assume that $\alpha $ and
$\gamma $ are algebra maps, while $\beta $ may be only a quasi-algebra map.

The following theorem provides the theoretical underpinnings for
noncommutative modeling of fiber squares.  It states that under certain reasonable
conditions, there exists a sort of push-out of $\iota $ and $\vp$ that is a model of the
pull-back $E\underset B\to{\times}X$.  

\proclaim{Theorem 0.3.2  \cite {DH1}}Given a commuting diagram over a field $\Bbbk$,
with squares as in diagrams (0.3.1) and (0.3.2),
$$\xymatrix{
(A,d)\odot (C,e)   \ar [d]^\beta _\simeq&(A,d)\ar[l]_(0.4)\iota\ar [r]^\vp\ar [d]^\alpha_\simeq&(\ba ,\bd)\ar [d]^\gamma_\simeq\\
C\sp *E& C\sp *B\ar [l]_{C^*q}\ar [r]^{C\sp *f}&C\sp *X}$$
in which $\bar A$ is a free algebra and $\vp $ admits a cochain
algebra section $\sigma $, there exist a twisted algebra extension 
$$\bar\iota : (\ba ,\bd)@>>>      (\ba ,\bd)\odot (C,e)$$
and a noncommutative model over $\Bbbk$ 
$$\delta :(\ba ,\bd )\odot (C,e)@>\simeq >> C\sp *(E\underset B\to{\times}X)$$ 
such that
\roster
\item $(\bar a\otimes 1)(1\otimes c)=(\vp \otimes 1)\bigl((\sigma(\bar
a)\otimes 1)(1\otimes c)\bigr)$ for all $\bar a$ in $\ba $ and $c$ in $C$;

\item if $D$ is the differential on $A\otimes C$, then the
differential $\bar D$ on $\ba\otimes C$ commutes with the right action of
$\ba $ and  is specified by
$\bar D (1\otimes c)=(\vp
\otimes 1)\bigl(D(1\otimes c)\bigr)$ for all $c$ in $C$;

\item $\delta (\bar a\otimes c)=C\sp *\bar f\circ \beta (\sigma (\bar
a)\otimes c)$ for all $\bar a$ in $\ba $ and $c$ in $C$.
\endroster   
\endproclaim

When we are not working over a field, as in this chapter, we cannot apply this theorem directly but have to employ more ad hoc methods, in order to obtain a result of this type.  In particular, defining the full algebra structure on $(\ba, \bd)\otimes (C,e)$ and then showing that $\delta$ is a quasi-algebra map can be delicate. 

\remark{Related work} As mentioned at the beginning of this subsection, our approach to algebraic modeling of fiber squares is analogous to the KS-extensions of rational homotopy theory, as developed by Sullivan \cite {FHT}.  The {\sl Adams-Hilton model}, which to any $1$-reduced  CW-complex $X$ associates a chain algebra $(AH(X),d)$ quasi-isomorphic to the cubical chains on $\Om X$ \cite {AH}, is another particularly useful tool for algebraic modeling.  The algebra $AH(X)$ is free on generators in one-to-one correspondence with the cells of $X$, and the differential $d$ encodes the attaching maps. 

In \cite {An}, Anick showed that the Adams-Hilton model could be endowed with a coproduct $\psi$, so that it became a {\sl Hopf algebra up to homotopy}.  He showed furthermore that  if $X$ is a finite 
$r$-connected CW complex of dimension at most $rp$, then there is a 
commutative cochain algebra $A(X)$ that is quasi-isomorphic to 
$C^*(X;\Bbb F_{p})$, the algebra of mod $p$ cochains on $X$. Using Anick's result, Menichi proved in \cite {Me} that if $i:X\hookrightarrow Y$ is an 
inclusion of finite $r$-connected CW complexes of dimension at most 
$rp$ and $F$ is the homotopy fiber of $i$, then the mod $p$ 
cohomology of $F$ is isomorphic as an algebra to $\operatorname 
{Tor}^{A(Y)}(A(X), \Bbb F_{p})$.
 
Other interesting algebraic models include the {\sl SHC-algebras} studied by Ndombol and Thomas \cite {NT} and {\sl $E_{\infty}$-algebras}, shown by Mandell to serve as models for $p$-complete homotopy theory \cite {Man}.  In particular, Mandell 
proved that the cochain functor $C^*(-;\overline{\Bbb F_{p}})$ embeds 
the category of nilpotent $p$-complete spaces onto a full subcategory 
of $E_{\infty}$-algebras.  He also characterized those 
$E_{\infty}$-algebras that are weakly equivalent to the cochains on 
a $p$-complete space.
\endremark


\head{1. Free loop spaces}\endhead

\noindent Consider the free loop fiber square for a simply-connected CW-complex of finite type,
$X$. 
$$\xymatrix{
\fls \ar[r]^j\ar@{->>}[d]^e&\cyl\ar @{->>}[d]^{(\ev _0,\ev_ 1)}\\
X \ar [r]^(0.4)\Delta & X\times X}$$
Here $\ev _t$ is defined by $\ev _t(\ell )=\ell (t)$, and
$\Delta $ is the diagonal. 

Our goal in this chapter is to construct canonically an associative cochain algebra $fls^*(X)$ together with a quasi-isomorphism $\xymatrix@1{fls ^*(X)\ar [r]^{\simeq}&C^*\Cal LX}$ that induces an isomorphism of algebras in cohomology. 

We construct the noncommutative model of $\fls $ over $\Bbb Z$ in
three steps.  First we find a model of $\Delta $, a relatively easy
exercise.  The second step, in which we define a twisted extension of cochain algebras that is a
 model of the topological fibration $(\ev _0,\ev_ 1)$, requires considerably more work.
Once we have obtained the models of $\Delta$ and $(\ev _0,\ev_ 1)$, we show that they can be ``twisted together," leading to a model for $\fls $. 

Since it is much easier to obtain precise, natural algebraic models for simplicial sets than for topological spaces, we begin this section by  constructing a useful, canonical simplicial model of the free loop space. We then apply the general theory of section 0.3 to building the desired algebraic free loop space model.

\subhead 1.1 A simplicial model for the free loop space\endsubhead
\subsubhead The general model\endsubsubhead
Let $\xymatrix@1{\Cal S_{\bullet}:\Cal Top\,\ar [r]&\,s\Cal Set}$ denote the singular simplicial set functor, which is a right adjoint to the geometric realization functor, with which it forms a Quillen equivalence.  Let $\xymatrix@1{\eta :Id \,\ar [r]&\,\ssimp |\cdot |}$ denote the unit of the adjunction. Recall that $\eta_L$ is always a weak equivalence \cite {May}.

Let $X$ be a $1$-connected space, and let $K$ be a $1$-reduced Kan complex such that $|K|$, the geometric realization of $K$, has the homotopy type of $X$.  For example, we could take $K=\Cal S_\bullet (X)$.  Let 
$$\xymatrix{K\ar [rr]^{\Delta}\ar@{ >->}[rd]_{\simeq}^i&&K\times K\\&P\ar @{->>} [ru]^p}$$ 
be a path object on $K$.  Set $\Cal L(K,P) := P\underset {K\times K}\to\times K$.  Let $\xymatrix@1{\bar p:\Cal L(K,P)\ar [r]&K}$ and $ \xymatrix @1{\overline\Delta  : \Cal L(K,P)\ar [r]&P}$ denote the canonical maps, i.e.,
$$\xymatrix{\Cal L(K,P)\ar [r]^{\overline\Delta }\ar [d]^{\bar p}&P\ar [d]^p\\
		   K\ar [r]^\Delta &K\times K}$$
is the pullback diagram.

\proclaim{Proposition 1.1.1} There is a weak equivalence $\xymatrix@1{\Cal L(K,P)\ar [r]&\, \ssimp (\Cal LX)}$.\endproclaim

\demo{Proof} 
Since $\ssimp$ is a right adjoint, it preserves limits.  It also preserves fibrations, as the right half of a Quillen equivalence. There is therefore a pullback diagram
$$\xymatrix{\ssimp(\Cal L|K|) \ar [rr] ^{\ssimp (j)}\ar @{->>}[d]^{\ssimp (e)}&&\ssimp (|K|^I)\ar @{->>} [d]^{(\ssimp (\ev_0),\ssimp (\ev_1))}\\ 
\ssimp (|K|)\ar [rr]^(.4){\ssimp (\Delta)}&&\ssimp (|K|)\times \ssimp (|K|)}$$ 
since $\ssimp (|K|\times |K|)\cong \ssimp (|K|)\times \ssimp (|K|)$.
Consider the following diagram, which commutes by naturality of $\eta$,
$$\xymatrix@R+=2pt{K\ar @{ >->}[d]^i_{\simeq}\ar [r]^(.4){\eta _K}_(.4){\simeq}\ar [dr]^(.6){\Delta}&\ssimp(|K|)\ar [rr] ^{\ssimp (s)}_{\simeq}\ar [drr]^(.6){\ssimp(\Delta)}&&\ssimp (|K|^I)\ar @{->>}[d]^{(\ssimp (\ev_0),\ssimp (\ev_1))}\\
P\ar@{->>} [r] ^(.4)p&K\times K\ar [rr]^(.4){\eta _{K^2}}_(.4){\simeq}&&\ssimp (|K|)\times \ssimp (|K|)}$$
where $\xymatrix@1{s:|K|\ar [r]&|K|^I}$ is the usual section, sending an element $x$ of $|K|$ to the constant path at $x$.  Since $i$ is an acyclic cofibration and $(\ssimp (\ev_0),\ssimp (\ev_1))$ is a fibration, there is a simplicial map $\xymatrix@1{\bar \eta:P\ar [r]&\ssimp (|K|^I)}$ such that
$$\xymatrix{K\ar @{ >->}[d]^i_{\simeq}\ar [rrr]^{\ssimp (s)\eta _K}_{\simeq}&&&\ssimp (|K|^I)\ar @{->>}[d]^{(\ssimp (\ev_0),\ssimp (\ev_1))}\\
P\ar [rrr] ^{\eta _{K^2}\circ p}\ar [urrr]^(.4){\bar\eta}&&&\ssimp (|K|)\times \ssimp (|K|)}$$
commutes.  Note that by ``2-out-of-3" $\bar\eta$ is a weak equivalence.

We have therefore a commutative diagram
$$\xymatrix{\Cal L(K,P)\ar [rdd]_{\eta _K\bar p}\ar [rrrd]^{\bar\eta\overline\Delta}\\
&\ssimp(\Cal L|K|) \ar [rr] ^(.4){\ssimp (j)}\ar @{->>}[d]^{\ssimp (e)}&&\ssimp (|K|^I)\ar @{->>} [d]^{(\ssimp (\ev_0),\ssimp (\ev_1))}\\ 
&\ssimp (|K|)\ar [rr]^{\ssimp (\Delta)}&&\ssimp (|K|)\times \ssimp (|K|)}$$ 
which implies that there is a simplicial map 
$$\xymatrix@1{\Cal L(K,P)\ar [r]^{\hat\eta}&\ssimp(\Cal L |K|)}$$
 such that $\ssimp (j)\hat \eta=\bar\eta \overline\Delta $ and $\ssimp (e)\hat\eta = \eta _K\bar p$. By the ``Cogluing Lemma" (cf. \cite {GJ: \S II.8}), it is clear that $\hat \eta $ is a weak equivalence.\qed 
\enddemo

For explicit computations to be possible, it is important to be able to build a simplicial model for $\Cal LX$ from a simplicial set that isn't necessarily Kan.  Suppose therefore that $K$ is any $1$-reduced simplicial set such that $|K|\simeq X$.  Let 
$$\xymatrix@1{K\,\cof ^j_{\simeq}&K' \fib & \{*\}}$$
be a fibrant replacement of $K$. 

Given path objects on $K$ and on $K'$,
$$\xymatrix@1{K\,\cof_{\simeq}^i&P\ar @{->>}[r]^(0.4)p&K\times K}$$
and
$$\xymatrix@1{K'\,\cof_{\simeq}^{i'}&P'\ar @{->>}[r]^(0.4){p'}&K'\times K',}$$
consider the commuting square
$$\xymatrix{K\ar @{ >->}[d]_{\simeq}^i\ar [r]_{\simeq}^j&K'\ar [r]_{\simeq}^{i'}&P'\ar @{->>}[d]^{p'}\\
P\ar [r]^p&K\times K\ar [r]^{j\times j}_{\simeq}&K'\times K'}$$
in which $j\times j$ is a weak equivalence.  Since $i$ is an acyclic cofibration and $p'$ is a fibration, there is a simplicial map $\xymatrix@1{\tilde \jmath:P\ar [r]&P'}$ such that 
$$\xymatrix{K\ar @{ >->}[d]_{\simeq}^i\ar [rr]_{\simeq}^{i'j}&&P'\ar @{->>}[d]^{p'}\\
P\ar [rr]^{(j\times j)p}\ar [urr]^(.4){\tilde \jmath}&&K'\times K'}$$
commutes.  By ``2-out-of-3", $\tilde \jmath$ is also weak equivalence.

If, as above, $\Cal L (K,P):= P\underset {K\times K}\to\times K$ and $\Cal L (K',P'):= P'\underset {K'\times K'}\to\times K'$, then by the universal property of pull-backs, there is a simplicial map $\xymatrix@1{\tilde \jmath':\Cal L (K,P)\ar [r]&\Cal L (K',P')}$ such that the following cube commutes.
$$\xymatrix{\Cal L(K,P)\ar [rr]\ar @{->>}[rd]^{\bar p}\ar [dd]^(.3){\tilde \jmath'}&&P\ar'[d]^(.6){\tilde \jmath}_(.6){\simeq} [dd]\ar@{->>}[dr]^p\\
&K\ar [dd]^(.3)j_(.3){\simeq}\ar [rr]&&K\times K\ar [dd]^(.3){j\times j}_(.3){\simeq}\\
\Cal L (K',P')\ar '[r][rr]\ar@{->>}[dr]^{\bar p'}&&P'\ar@{->>} [dr]^{p'}\\
&K'\ar [rr]&&K'\times K'}$$
As in the proof of Proposition 1.1.1, the ``Cogluing Lemma" implies that $\tilde \jmath'$ is a weak equivalence.  In particular, $\H^*(\Cal L(K,P))\cong \H^*(\Cal L (K',P'))$, as algebras.  We have thus established the following result.

\proclaim{Theorem 1.1.2} Let $X$ be a $1$-connected space, and let $\Cal LX$ be the free loop space on $X$. Let $K$ be any $1$-reduced simplicial set such that $|K|$ has the homotopy type of $X$.  Let $\xymatrix@1{K\,\cof_{\simeq}^i&P\ar @{->>}[r]^p&K\times K}$ be a path object on $K$, and set $\Cal L (K,P):= P\underset {K\times K}\to\times K$.  Then there is simplicial weak equivalence $$\xymatrix@1{\Cal L (K,P)\ar [r]^{\simeq}&\,\ssimp (\Cal LX)}.$$  In particular, $\H^*(\Cal L(K,P))\cong \H^*(\Cal LX)$ as algebras.\endproclaim

\subsubhead Choosing the free loop model functorially\endsubsubhead
\noindent Let $K$ be a $1$-reduced  simplicial set.  Let $GK$ denote the Kan loop group on $K$, as defined in the preface. Let $PK$ denote the twisted cartesian product with structure group $GK\times GK$
$$PK:=GK\underset \tau\to \times (K\times K),$$
where $\xymatrix@1{\tau :K\times K\ar [r]&\,GK\times GK:(x,y)\,\ar @{|->}[r]&\,(\bar x,\bar y)}$ and 
$GK\times GK$ acts on $GK$ by $(v,w)\cdot u:=vuw^{-1}$. It is easy to verify that $\tau$ satisfies the conditions of a twisted cartesian product. Observe that, in particular,
$$\del _0(u, (x,y))=(\bar x\cdot \del _0 u\cdot \bar y^{-1}, (\del _0x, \del _0y)).$$

\proclaim{Proposition 1.1.3} Let $\xymatrix@1{p:PK\ar[r]&K\times K}$ denote the projection map, and let 
$$\xymatrix@1{i:K\ar [r]&PK:x\,\ar @{|->}[r]&\,(e,x,x)}.$$
Then $\xymatrix@1{K\ar [r]^i&PK\ar [r]^p&K\times K}$ is a path object on $K$.\endproclaim

\demo{Proof} Since $PK$ is a twisted cartesian product and $GK$ is a Kan complex, $p$ is a Kan fibration.  Furthermore, $i$ is a simplicial cofibration, since it is obviously injective.  We need therefore only to show that $i$ is a weak equivalence. 

If $K$ is a Kan complex, then there is a long exact sequence of homotopy groups
$$\xymatrix@1{\cdots\ar [r]&\pi _{n+1} (K\times K)\ar [r]^(0.6){\delta_{\tau}}&\pi _n(GK)\ar [r]&\pi_n(PK)\ar [r]&\pi _n(K\times K)\ar [r]&\cdots},$$
where $\delta_{\tau} $ is the connecting homomorphism, defined by 
$\delta_{\tau} \bigl([x,y]\bigr)=[\bar x \cdot \bar y^{-1}]$ (cf., \cite{GJ:\S I.7}).  Comparing this long exact sequence with that obtained from the universal acyclic twisted cartesian product $GK\underset \eta \to \times K=EK$
$$\xymatrix@1{\cdots\ar [r]&\pi_{n+1}(EK)\ar [r] &\pi _{n+1} (K)\ar [r]^{\delta_{\eta}}&\pi _n(GK)\ar [r]&\pi_n(EK)\ar [r]&\cdots,}$$
we obtain that $\delta _{\eta}$ is an isomorphism given by $\delta _{\eta}\bigl ([x]\bigr )=[\bar x]$, and therefore that $\delta _{\tau}$ is surjective. Furthermore, $\ker \delta _{\tau}= \text{Im}\, \pi _{n+1}\Delta$, since
$$\xymatrix@1{[x,y]\in \ker \delta _\tau\ar @{<=>}[r]&[\bar x\cdot \bar y^{-1}]=[e]\ar @{<=>}[r]&[\bar x]=[\bar y]\ar @{<=>}[r]&[x]=[y]\ar @{<=>}[r]&[x,y]=[x,x]}$$
where the second equivalence comes from multiplying both sides of the equation by $[\bar y]$, using the multiplication induced by that on $GK$, while the third equivalence is due to the fact that $\delta _{\eta}$ is an isomorphism.  

The first long exact sequence therefore breaks up into short exact sequences
$$\xymatrix@1{0\ar [r]&\pi _{n+1}(PK)\ar [rr]^(0.45){\pi _{n+1}p}&&\pi _{n+1} (K\times K)\ar [r]^(0.6){\delta_{\tau}}&\pi _n(GK)\ar [r]&0},$$
which are split, since $\delta _{\tau}$ has an obvious section $\sigma$ defined by $\sigma \bigl([\bar x])=[x,*]$, where $*$ is (an iterated degeneracy of) the basepoint of $K$.  Consequently $\pi _{n+1}p$ has a left inverse $\rho$, and so
$$\pi _{n+1} i=\rho \pi _{n+1}p\pi _{n+1}i=\rho \pi _{n+1} \Delta.$$
Since $\im \pi _{n+1} \Delta=\ker \delta _\tau=\im \pi _{n+1} p$ and the restriction of $\rho$ to $\im \pi _{n+1} p$ is surjective, $\pi _{n+1} i$ is surjective as well. Furthermore $\pi _{n+1} i$ is necessarily injective, as $\pi _{n+1} pr_1\pi _{n+1}p$ is a left inverse to $\pi _{n+1} i$, where $\xymatrix@1{pr_1:K\times K\ar [r]&K}$ is the projection map onto the first factor.  Thus $\pi _{n+1} i$ is an isomorphism for all $n$, i.e., $i$ is a weak equivalence.

If $K$ is not a Kan complex, consider a fibrant replacement $\xymatrix@1{K\,\cof ^j_{\simeq}&K' \fib & \{*\}}$, which, by the naturality of the $P$-construction, induces a map of twisted cartesian products
$$\xymatrix{GK\ar @{ >->}[d]_{\simeq}^{Gj}\ar [r] &PK\ar [d] ^{Pj}\ar @{->>}[r]^p&K\ar [d]_{\simeq }^j\\
GK'\ar [r]&PK'\ar @{->>}[r]&K'}$$
in which $Gj$ is an acyclic cofibration, since $G$ is the left member of a Quillen equivalence.  A 5-Lemma argument applied to the long exact sequences in homotopy of the realizations of these fibrations then shows that $Pj$ is a weak equivalence.  Again by naturality, the square 
$$\xymatrix{K\ar [d]^j_{\simeq}\ar [r] ^i&PK\ar [d]_{\simeq}^{Pj}\\
K'\ar [r]^{i'}_{\simeq}&PK'}$$
commutes and so, by ``2-out-ot-3", $i$ is also a weak equivalence.
\qed\enddemo 

We now use this functorial path object construction to define a functorial, simplicial free loop space model.

\definition{Definition} Given a $1$-reduced simplicial set $K$, the {\sl canonical free loop construction} on $K$ is the the simplicial set
$$\Cal L K:=\Cal L(K,PK)=PK\underset K\times K\to \times K.$$\enddefinition

Observe that $\Cal LK$ is the twisted cartesian product $GK\underset \bar\tau\to \times K$, with structure group $GK\times GK$, where $GK\times GK$ acts on $GK$ as before and
$$\xymatrix@1{\bar\tau:K\ar [r]&GK\times GK:x\ar @{|->}[r]&(\bar x, \bar x)}$$
so that 
$$\del _0(w,x)=(\bar x\cdot \del _0w\cdot \bar x^{-1}, \del _0x).$$

According to Theorem 1.1.2, there is a weak equivalence $\xymatrix@1{\Cal LK\ar [r]^(.4){\simeq}&\,\ssimp (\Cal L|K|)}$ and therefore an algebra isomorphism $\H^*(\Cal LK)\cong \H^*(\Cal L |K|)$.

\subhead 1.2 The multiplicative free loop space model\endsubhead
\noindent In this section we apply the methods described in section 0.3 to constructing naturally a noncommutative model for the canonical free loop construction $\Cal L K$ on a $1$-reduced simplicial set of finite-type $K$.  By Theorem 1.1.2 we obtain therefore  a noncommutative model for the free loop space $\Cal LX$ on a $1$-connected space $X$ with the homotopy type of a finite-type CW-complex.  

We begin by constructing specific explicit models of the diagonal map$$\xymatrix@1{\Delta :K\ar[r]&\,K\times K}$$ and of the path fibration $$\xymatrix@1{p:PK\ar [r]&\,K\times K,}$$ which we then twist together appropriately, in order to obtain a model of the free loop construction.

\subsubhead The diagonal map\endsubsubhead
\noindent Let $\xymatrix@1{\varepsilon :\Om\bc\ar [r]&Id}$ denote the counit of the cobar-bar adjunction, which is a quasi-isomorphism for each $1$-connected cochain algebra $(A,d)$. Let $$\xymatrix@1{\gamma=\varepsilon_K:\Om\bc C^*K\ar [r]^(0.7){\simeq}&C^*K}.$$
Let $\alpha$ denote the following composition.
$$\xymatrix{\Om(\bc C^*K\otimes \bc C^*K)\ar [r]^{\Om (\rho^{\sharp})}\ar [drrr]_\alpha&\Om\bc (C^*K\otimes C^*K)\ar [rr]^{\Om (\tom f_{K,K})^{\sharp}}&&\Om \bc C^*(K\times K)\ar [d]^{\varepsilon _{K\times K}}_{\simeq}\\&&&C^*(K\times K)}$$
Here $\rho ^{\sharp}$ is the dual of the map defined in equation (0.2.1), while $f_{K,K}$ is again the Alexander-Whitney map, as in section 0.2.  Note that we are relying on the fact that $K$ is of finite-type, in writing, e.g., $\bc C^*K$ for the dual of $\Om C_*K$.

Given these definitions, it is an easy exercise, using the naturality of $\varepsilon$, to show that 
$$\xymatrix{\Om(\bc C^*K\otimes \bc C^*K)\ar [d]_{\simeq}^{\alpha}\ar [rr]^(0.6){\Om \psi _K^{\sharp}}&&\Om \bc C^*K\ar [d] _{\simeq}^{\gamma}\\
C^*(K\times K)\ar [rr]^{C^*\Delta}&&C^*K}\tag 1.2.1$$
commutes and is therefore a model of $\Delta$. Observe that if
$$\xymatrix{\iota _k: \Om\bc C^*K \ar[r]&\Om(\bc C^*K\otimes \bc C^*K)}$$
is the inclusion of the $k$th tensor factor ($k=1,2$), then $C^*pr_k\circ\gamma =\alpha\iota _k$.

\subsubhead The path fibration\endsubsubhead
\noindent Henceforth, we employ the following useful notation. 

\definition{Notation} For any $a,b\in \bc C^*K$, let $a\star b:=\psi _K^{\sharp}(a\otimes b)$. Furthermore, let $\vp=\Om\psi _K^{\sharp}$. \enddefinition

We begin by defining a certain twisted extension of $\Om(\bc C^*K\otimes \bc C^*K)$ by $\bc C^*K$.  Since the former is a model of $K\times K$, while the latter is a model of $GK$, it is reasonable to expect to be able to construct a model of $P K$ of this form.  Once we have the explicit definition of the twisted extension, we show that it is indeed a model of $PK$.

\definition{Definition of $\Om(\bc C^*K\otimes \bc C^*K)\odot \bc C^*K$}

\noindent Let $\Om(\bc C^*K\otimes \bc C^*K)\odot \bc C^*K$ be the twisted algebra extension determined by the following conditions.
\roster

\item The left action of $\Om(\bc C^*K\otimes \bc C^*K)$  is specified recursively by
$$\align
\bigl (\si (a\otimes b)&\otimes  1\bigr)\cdot (1\otimes sx_1|\cdots |sx_n)\\=&\si (a\otimes b)\otimes sx_1|\cdots |sx_n\\
&+(-1)^{b\cdot \theta _n}\biggl(\si \bigl(a\otimes (sx_1|\cdots |sx_n\star b)\bigr)-\si\bigl ((a\star sx_1|\cdots |sx_n)\otimes b\bigr)\biggr)\\
&+\sum _{j=1}^{n-1}\biggl[(-1)^{b\cdot \theta _j}\si \bigl((a\star sx_1|\cdots |sx_j)\otimes b\bigr)\cdot (1\otimes sx_{j+1}|\cdots |sx_n)\\
&\qquad +(-1)^{(b+\theta _j)(\theta_n-\theta _j)}\si \bigl(a\otimes (sx_{j+1}|\cdots |sx_n)\star b\bigr)\otimes sx_1|\cdots |sx_j\biggr]\endalign$$
where $\theta_j=j+\sum _{i\leq j}|x_i|$.
\item The restriction of the differential $D$ of the twisted extension to $1\otimes \bc C^*K$ is specified recursively by 
$$\align
D(1\otimes sx_1|\cdots |sx_n)=&1\otimes d_{\bc}(sx_1|\cdots |sx_n)\\
 &+\si (sx_1|\cdots |sx_n\otimes 1)-\si (1\otimes sx_1|\cdots |sx_n)\\
 &+\sum _{j=1}^{n-1}\bigl[\si (sx_1|\cdots |sx_j\otimes 1)\cdot (1\otimes sx_{j+1}|\cdots |sx_n) \\
 &\qquad -(-1)^{\theta _j(\theta _n-\theta j)}\si (1\otimes sx_{j+1}|\cdots |sx_n)\otimes sx_1|\cdots |sx_j\bigr]\endalign$$
 with $\theta _j$ as above.
 
 Observe that according to this definition, if $c=sx_1|\cdots |sx_n$, then $D(1\otimes c)$ is of the following form.
 $$D(1\otimes c)=1\otimes d_{\bc}c+\sum _{i<j} \si \bigl(\lambda _{ij,k}(c)\otimes \lambda _{ij}^k(c)\bigr)\otimes sx_i|\cdots |sx_j$$
 where we are applying the Einstein summation convention (cf., section 0.1).
 \item Let $a=sx_1|\cdots |sx_m,b=sy_1|\cdots |sy_n\in \bc C^*K$. Then we define 
 $$\align
 (1\otimes a)&(1\otimes b)\\
 =&1\otimes a\star b\\
 &+\sum \Sb i<j\\ p<q\endSb\biggl[\pm\si\bigl(\Theta^{\sharp} (\lambda _{ij,k}(a)\otimes \lambda _{pq,r}(b))\otimes \lambda _{ij}^k(a)\star \lambda _{pq}^r(b)\bigr)\\
&\qquad\quad\;\pm \si (\lambda _{ij,k}(a)\star \lambda _{pq,r}(b)\otimes \Theta^{\sharp} (\lambda _{ij}^k(a)\otimes \lambda _{pq}^r(b))\bigr)\\
&\qquad\quad\;\pm \si \bigl(\Theta^{\sharp} (\lambda _{ij,k}(a)_s\otimes \lambda _{pq,r}(b)_t)\otimes \lambda _{ij}^k(a)_u\star \lambda _{pq}^r(b)_v\bigr)\cdot\\
&\qquad\qquad\qquad\qquad  \cdot\si (\lambda _{ij,k}(a)^s\star \lambda _{pq,r}(b)^t\otimes \Theta^{\sharp} (\lambda _{ij}^k(a)^u\otimes \lambda _{pq}^r(b)^v)\bigr)\biggr]\\
&\qquad\qquad\qquad\qquad\qquad\qquad\qquad\qquad\qquad\otimes (sx_i|\cdots |sx_j)\star( sy_p|\cdots |sy_q) 
\endalign$$
where we have suppressed the relatively obvious, though horrible to specify, signs given by the Koszul convention (cf., section 0.1). Here 
$\Theta^{\sharp}$ is the coderivation homotopy for the commutativity of the multiplication in $\bc C^*K$, dual of the derivation homotopy $ \Theta$ for the cocommutativity of $\psi _K$, and the extra lower and upper indices in the last summand indicate factors of the coproduct evaluated on the terms in question, e.g.,
$\Delta (\lambda _{ij,k}(a))= \lambda _{ij,k}(a)_s\otimes \lambda _{ij,k}(a)^s$.

Observe that, in particular, 
$$(1\otimes a)(1\otimes b) \in T^{\leq 2}\si \bc C^*K\otimes \bc C^*K.$$
\endroster
\enddefinition

Dupont and Hess showed in \cite {DH4} that the conditions above specify an associative cochain algebra, when the product $\psi_K^{\sharp}$ is strictly commutative, so that condition (3) reduces to $(1\otimes a)\cdot(1\otimes b)=1\otimes (a\star b)$.  The general, homotopy-commutative case was established by Blanc in his thesis \cite {Bl}.

Extend $\xymatrix@1{\vp :\Om (\bc C^*K\otimes \bc C^*K)\ar[r]&\Om\bc C^*K}$ to 
$$\xymatrix@1{\overline\vp:\Om (\bc C^*K\otimes \bc C^*K)\odot \bc C^*K\ar[r]&\Om\bc C^*K}$$
by $\overline\vp (1\otimes w)=0$ for all $w\in (\bc C^*K)^+$.  It is easy to see that $\overline \vp$ is a differential map and that, as explained in \cite {DH2}, a straightforward spectral sequence argument shows as well that $\overline\vp$ is a quasi-isomorphism.   If  $\psi_K^{\sharp}$ is strictly commutative, then $\overline\vp$ is a strict algebra map, as proved in \cite {DH4}. More generally, Blanc showed in \cite {Bl} that $\overline\vp$ is a quasi-algebra map when $\psi_K^{\sharp}$ is commutative up to a derivation homotopy. We have thus established the following result.

\proclaim{Lemma 1.2.1} There is a commuting diagram 
$$\xymatrix{\Om (\bc C^*K\otimes \bc C^*K)\ar[rr]^{\vp}\ar [dr]^{\iota}&&\Om\bc C^*K\\
&\Om (\bc C^*K\otimes \bc C^*K)\odot \bc C^*K\ar [ur]^{\overline\vp}_{\simeq}}$$
where $\iota $ is the inclusion, $\overline\vp (1\odot w)=0$ for all $w\in (\bc C^*K)^+$, and $\overline\vp $ is a quasi-algebra quasi-isomorphism.\endproclaim

Now consider the following commutative diagram
$$\xymatrix{\Om(\bc C^*K\otimes \bc C^*K)\ar [r]^{\alpha}_{\simeq}\ar [d]^{\iota}&C^*(K\times K)\ar [r]^{C^*p}&C^*PK\ar [d]^{C^*i}\\
\Om (\bc C^*K\otimes \bc C^*K)\odot \bc C^*K\ar [r]^(0.7){\overline\vp}_(0.7){\simeq}&\Om\bc C^*K\ar [r]^{\gamma}_{\simeq}&C^*K}$$
which satisfies the conditions of Proposition 0.3.1.  Hence, we can lift $\overline\vp \gamma$ through $C^*i$, obtaining a quasi-algebra map
$$\xymatrix@1{\beta :\Om(\bc C^*K\otimes C^*K)\odot \bc C^*K\ar [r]&C^*PK}$$
such that $\beta\iota=C^*p\circ \alpha$ and $C^*i\beta=\gamma\circ\overline\vp$. We can therefore take
$$\xymatrix{\Om(\bc C^*K\otimes \bc C^*K)\ar [d]^{\alpha}_{\simeq}\ar [r]^(0.4){\iota}&\Om (\bc C^*K\otimes \bc C^*K)\odot \bc C^*K\ar [d]_{\simeq}^{\beta}\\
C^*(K\times K)\ar [r]^{C^*p}&C^*PK}\tag 1.2.2$$
as an algebraic model of the path fibration.

\subsubhead The free loop space model\endsubsubhead
\noindent We now twist together the models (1.2.1) and (1.2.2), in the spirit of Theorem 0.3.2.

\proclaim{Theorem 1.2.2 \cite {DH4},\cite {Bl} }There is  a twisted algebra extension 
$$\xymatrix{\bar\iota:\Om\bc C^*K\ar [r]&\Om\bc C^*K\odot \bc C^*K}$$ and a quasi-algebra quasi-isomorphism
$$\xymatrix{\delta:\Om \bc C^*K\odot \bc C^*K\ar [r]^(0.6){\simeq}& C^*(\Cal L K)}$$ 
defined as follows.
\roster
\item The left action of $\Om\bc C^*K$  is specified recursively by
$$\align
\bigl (\si a&\otimes  1\bigr)\cdot (1\otimes sx_1|\cdots |sx_n)\\=&\si a\otimes sx_1|\cdots |sx_n\\
&-\sum _{j=1}^{n-1}\biggl[\si \bigl(a\star (sx_1|\cdots |sx_j))\cdot (1\otimes sx_{j+1}|\cdots |sx_n)\\
&\qquad +(-1)^{(\theta _j)(\theta_n-\theta _j)}\si \bigl(a\star (sx_{j+1}|\cdots |sx_n)\bigr)\otimes sx_1|\cdots |sx_j\biggr]\endalign$$
where $\theta_j=j+\sum _{i\leq j}|x_i|$.
\item The restriction of the differential $\overline D$ of the twisted extension to $1\otimes \bc C^*K$ is specified recursively by 
$$\align
\overline D(1\otimes sx_1|\cdots |sx_n)=&1\otimes d_{\bc}(sx_1|\cdots |sx_n)\\
 &+\sum _{j=1}^{n-1}\bigl[(\si (sx_1|\cdots |sx_j)\otimes 1)\cdot (1\otimes sx_{j+1}|\cdots |sx_n) \\
 &\qquad -(-1)^{\theta _j(\theta _n-\theta j)}\si (sx_{j+1}|\cdots |sx_n)\otimes sx_1|\cdots |sx_j\bigr]\endalign$$
 with $\theta _j$ as above.
 \item Let $a=sx_1|\cdots |sx_m,b=sy_1|\cdots |sy_n\in \bc C^*K$. Then we define 
 $$\align
 (1\otimes a)&(1\otimes b)\\
 =&1\otimes a\star b\\
 &+\sum \Sb i<j\\ p<q\endSb (-1)^{\zeta _{ij,pq}}\si\biggl(\Theta ^{\sharp} \bigl(\lambda _{ij,k}(a)\star \lambda _{ij}^k(a)\otimes \lambda _{pq,r}(b)\bigr)\star \lambda _{pq}^r (b)\biggr)\\
 &\qquad\qquad\qquad\qquad\qquad\qquad\qquad\qquad\otimes(sx_i|\cdots |sx_j)\star( sy_p|\cdots |sy_q). 
 \endalign$$
 Here we are using the same notation as in the definition of the path space model, and 
 $$\zeta_{ij,pq}=\bigl( |\lambda _{pq,r}(b)|+|\lambda _{pq}^r(b)|\big)\cdot (\theta _j-\theta _{i-1}).$$
 \item For all $a,b\in \bc C^*K$,
 $$\align
 \delta (\si a \otimes b):&=C^*\overline\Delta \circ \beta (\si (a\otimes 1)\otimes b)\\
 &=C^*\overline\Delta\circ \beta (1\otimes b)\cdot C^*\bar p\circ\gamma (\si a).
 \endalign$$
 \endroster\endproclaim  
 This theorem was proved in the strictly commutative case in \cite {DH4} and in the homotopy-commutative case in \cite {Bl}.

In the article \cite {BH} based on part of Blanc's thesis, we use this model to compute the free loop space cohomology algebra of a space $X$ not in the ``Anick" range, i.e., where the space does not have a strictly commutative model.  This is the first time an explicit calculation of this type has been carried out.

\remark{Related work}  Our free loop space model construction is patterned on the construction for rational spaces, due to Sullivan and Vigu\'e \cite {SV}.  They used their construction to prove the rational {\sl Closed Geodesic Conjecture}.

Kuribayashi, alone \cite {K} and with Yamaguchi \cite {KY}, has applied the Eilenberg-Moore spectral sequence to calculating the cohomology algebra of certain free loop spaces.

Let $\frak C_*(A)$ denote the Hochschild complex of a (co)chain algebra $A$. Let $N^*X$ denote the normalized singular cochains on a topological space $X$ with coefficients in a field $\Bbbk$. Ndombol and Thomas  have applied SHC-algebra methods to showing that there is a natural map of cochain complexes $\xymatrix@1{\frak C_*(N^*X) \ar [r]&C^*(\Cal LX)}$ inducing an isomorphism of graded algebras in cohomology.  In a similar vein, Menichi proved in \cite {Me2} that $\H^*(\Cal LX)$ is isomorphic to the Hochschild cohomology of the singular chains on $\Omega X$.

Very recently B\"okstedt and Ottosen have developed yet another promising method for calculating free loop space cohomology, via a Bousfield spectral sequence \cite {BO1}.  Their method, which they have applied explicitly to spaces $X$ with $\H^*(X;\Bbb F_2)$ a truncated polynomial algebra on one generator, allows them to obtain the module structure of $\H^*(X;\Bbb F_2)$ over the Steenrod algebra.
\endremark

\subhead 1.3 The free loop model for topological spaces\endsubhead

\noindent In this section we adapt the free loop model above to the case of topological spaces, in order to facilitate construction of a model of the $S^1$-homotopy orbits on the free loop space. We sacrifice perhaps a bit of the multiplicative structure of the model in section 1.2, but since we are interested in this chapter only in the linear structure of the homotopy orbit cohomology, this is not a great loss.

 As usual let $X$ be a $1$-connected space with the homotopy type of a finite-type CW-complex, and let $K$ be a finite-type, $1$-reduced simplicial set such that $|K|\simeq X$. We again adopt the notation $\vp$ for $\Om \psi _K^{\sharp}$, where $(\Om C_*K, \psi_K)$ is the canonical Adams-Hilton model of section 0.2. Recall furthermore that $\varepsilon_K$ denotes the unit of the cobar-bar adjunction.

The unit $\eta $ of the $(|\cdot|, \ssimp)$-adjunction induces an injective quasi-isomorphism of chain coalgebras
 $$\xymatrix{C_*\eta _K:C_*K\ar [r]^(0.6)\simeq&S_*|K|}$$
 that dualizes to a surjective quasi-isomorphism of cochain algebras
 $$\xymatrix{C^*\eta _K:S^*|K|\ar [r]^(0.6){\simeq}&C^*K}.$$
 
 Furthermore,  for all spaces $Y$, subdivision of cubes defines a quasi-isomorphism of chain coalgebras
 $$\xymatrix{CU_*Y\ar [r]^\simeq&S_*Y,}$$
 which, upon dualization, gives rise to a quasi-isomorphism of cochain algebras
 $$\xymatrix{S^*Y\ar [r]^\simeq&CU^*Y.}$$ 

We begin the topological adaptation of the model for the simplicial free loop space by carefully specifying a model of the diagonal map, directly in the cubical cochains, which is where we know how to work explicitly with circle actions, as we show in section  2.

\proclaim{Theorem 1.3.1}  There is a noncommutative model of the diagonal map
$$\xymatrix{\Om(\bc C^*K\otimes \bc C^*K) \ar [d]^{\tilde\alpha}_\simeq\ar [r]^(0.6){\vp}&\Om \bc C^*K\ar [d]^{\tilde \gamma}_\simeq\\
CU^*(X\times X)\ar[r]^{CU^*\Delta}&CU^*X}$$
such that 
\roster
\item $\tilde \alpha\circ\iota _k=C^* pr_k\circ\tilde\gamma$ for $k=1,2$, where $$\xymatrix{\iota _k: \Om\bc C^*K \ar[r]&\Om(\bc C^*K\otimes \bc C^*K)}$$ is the inclusion into the $k$th tensor factor and 
$pr_k$ is the projection onto the  $k$th factor, and 
\item $\ker \varepsilon _K\subseteq \ker \tilde \gamma$.\endroster\endproclaim

We do not claim here that $\tilde\gamma$ and $\tilde\alpha$ are strict algebra maps.  In fact, it is probably impossible in general for $\tilde\gamma$ simultaneously to be an algebra map and to satisfy $\ker \varepsilon _K\subseteq \ker \tilde \gamma$. We show in the proof, however, that $\tilde\gamma$ and $\tilde\alpha$ are at least quasi-algebra maps.

\demo{Proof} Consider the following commutative diagram, in which the solid arrows are given and are all cochain algebra maps. We will explain step-by-step the construction and properties of the dashed arrows.
$$\xymatrix{
(\Om\bc C^*K)^{\oplus 2}\ar @{ >->}[dd]_{\iota =\iota _1+\iota _2}\ar @{-->}@/^1.25pc/[drr]^{(ii)\quad \exists \alpha'}\\
&&S^*|K|^2\ar @{->>}[dl]_{C^*\eta _{K^2}}^\simeq\ar @{->>}[dd]^{S^*\Delta}\\
\Om(\bc C^*K^{\otimes 2})\ar [dd]_{\vp}\ar [r]^(0.6)\alpha_(0.6)\simeq \ar @{-->}@/^1.5pc/[urr]^{(iii)\quad \exists\alpha ''\simeq}&C^*K^2\ar @{->>}[dd]^{C^*\Delta} \\
&&S^*|K|\ar @{->>}[dl]_(0.6){C^*\eta _K}^(0.6)\simeq \\
\Om\bc C^*K\ar [r]^\gamma_\simeq \ar @{-->}@/_4.5pc/[urr]_{(i)\quad\exists\hat\gamma\simeq}&C^*K}$$

\noindent {\sl Step (i)} Consider the following commutative diagram of cochain complexes.
$$\xymatrix{\ker \ve _K\ar @{ >->}[d]^{\text{incl.}}\ar [rr]^0&&S^*|K|\ar @{->>}[d]^{C^*\eta _K}_\simeq\\
\Om\bc C^*K\ar [rr]^{\gamma =\ve_K}_{\simeq}&&C^*K}$$
Since $C^*\eta _K$ is a surjective quasi-isomorphism and the inclusion map on the left is a free extension of cochain complexes, there exists $\xymatrix@1{\hat\gamma : \Om\bc C^*K\ar [r]&S^*|K|}$ such that $\ker\ve _K\subset \ker\hat\gamma$ and $S^*|K|\circ \hat\gamma =\ve _K$.  In particular $\hat\gamma$ is a quasi-isomorphism by ``2-out-of-3" and $\H^*\hat\gamma=(\H^*\eta _K)^{-1}\circ \H^*\gamma$ is a map of algebras, i.e., $\hat\gamma$ is a quasi-algebra map.
\smallskip
\noindent {\sl Step (ii)} Let $\alpha '=S^*pr_1\circ \hat\gamma +S^*pr_2\circ \hat\gamma$. Since $\vp\iota _k=Id$ for $k=1,2$ and $S^*\Delta\circ S^*pr_k=Id$, the diagram
$$\xymatrix{(\Om\bc C^*K)^{\oplus 2}\ar [r]^{\alpha '}\ar [d]^{\vp\iota}&S^*|K|^2\ar [d]^{S^*\Delta}\\
\Om\bc C^*K\ar [r]^{\hat\gamma}&S^*|K|}$$
commutes, by the universal property of the direct sum.
\smallskip
\noindent{\sl Step (iii)} Consider finally the commutative diagram below of cochain complexes and maps.
$$\xymatrix{(\Om\bc C^*K)^{\oplus 2}\ar [r]^(0.6){\alpha '}\ar [d]^{\iota}&S^*|K|^2\ar @{->>}[d]^{C^*\eta _{K^2}}_\simeq\\
\Om(\bc C^*K^{\otimes 2})\ar [r]^(0.6){\alpha}&C^*K^2}$$
As usual, since $C^*\eta _{K^2}$ is a surjective quasi-isomorphism and $\iota$ is a free extension of cochain complexes, $\alpha$ lifts to $\xymatrix@1{\alpha'':\Om(\bc C^*K^{\otimes 2})\ar [r]&S^*|K|^2}$, which is a quasi-algebra quasi-isomorphism such that $\alpha ''\iota =\alpha '$ and $C^*\eta _{K^2}\circ \alpha ''=\alpha$.

Notice that 
$$\align C^*\eta _K\circ S^*\Delta \circ \alpha ''=&C^*\Delta\circ  C^*\eta _{K^2}\circ \alpha ''\\
=&C^*\Delta\circ \alpha\\
=&\gamma\circ\vp\\
=&C^*\eta _K\circ \hat\gamma\circ\vp.
\endalign$$
We would like to have $S^*\Delta \circ \alpha ''= \hat\gamma\circ\vp$, but, as we show next, all we can be sure of is that $S^*\Delta \circ \alpha ''\simeq \hat\gamma\circ\vp$, as cochain maps.  We then apply the homotopy extension property of free extensions to ``fix" $\alpha ''$ and thus conclude the proof.

For any cochain complex, $(A,d)$, with free underlying graded abelian group,  let $I(A,d)=(A\oplus A\oplus sA, D)$ denote its canonical cylinder.   Let $\xymatrix@1{\frak j_k:A\ar [r]& I(A,d)}$ denote the inclusion onto the $k$th summand ($k=1,2$).  Note that $\frak j_k$ is always a quasi-isomorphism. 

Let $(P,d)$ be the cochain complex such that the following diagram is a push-out in the category of cochain complexes.
$$\xymatrix{\Om(\bc C^*K)^{\oplus 2}\oplus\Om(\bc C^*K)^{\oplus 2}\ar @{>->}[d]^{\iota\oplus\iota}\ar [rr]^(.6){\frak j_1+\frak j_2}&& I(\Om(\bc C^*K)^{\oplus 2})\ar [d] \\
\Om(\bc C^*K^{\otimes 2})\oplus\Om(\bc C^*K^{\otimes 2})\ar [rr]&&(P,d)}$$
We then have a commuting diagram
$$\xymatrix{(P,d)\ar [rrr] ^{(S^*\Delta\circ\alpha''+\hat\gamma\vp)+G}\ar  [d]^{\text{incl.}}&&&S^*|K|\ar @{->>}[d]_\simeq^{C^*\eta_K}\\
I\bigl(\Om(\bc C^*K^{\otimes 2})\bigr)\ar [rrr]^H&&&C^*K}$$
where $G$ is the constant homotopy (i.e., $G(s\Om(\bc C^*K)^{\oplus 2})=\{0\}$) from $S^*\Delta \circ \alpha '=\hat\gamma\vp\iota$ to itself and $H$ is the constant homotopy from $C^*\Delta\circ\alpha=\gamma\vp$ to itself. Again, since $C^*\eta _K$ is a surjective quasi-isomorphism and the inclusion on the left is a free extension, we can lift the homotopy $H$ to 
$$\xymatrix{H':I\bigl(\Om(\bc C^*K^{\otimes 2})\bigr)\ar [r]&S^*|K|}$$ 
such that $$H'|_{(P,d)}=(S^*\Delta\circ\alpha''+\hat\gamma\vp)+G,$$ i.e., $H'$ is a homotopy from $S^*\Delta\circ\alpha''$ to $\hat\gamma\vp$ that is constant on the subcomplex $\Om(\bc C^*K)^{\oplus 2}$.

Let $\xymatrix@1{G':I(\Om(\bc C^*K)^{\oplus 2})\ar [r]&S^*|K|^2}$ be the constant homotopy from $\alpha '$ to itself.  In particular, $S^*\Delta\circ G'=G$.  Consider the push-out diagram of cochain complexes
$$\xymatrix{\Om(\bc C^*K)^{\oplus 2}\ar @{>->} [d]^{\iota}\cof ^{\frak j_1}_\simeq &I(\Om(\bc C^*K)^{\oplus 2})\ar @{>->} [d]\ar @/^1.5pc/ [rdd]^{G'}\\
\Om(\bc C^*K^{\otimes 2})\cof ^{\hat\frak j_1}_{\simeq}\ar @/_1.25pc/[rrd]_{\alpha ''}^\simeq &(Q,d)\ar @{-->}[dr]^{\alpha ''+G'}_\simeq\\
&&S^*|K|^2}$$
in which all arrows are free extensions, since the push-out of a cofibration is a cofibration. Furthermore,  $\hat\frak j_1$ is a quasi-isomorphism because the push-out of an acyclic cofibration is an acyclic cofibration. 
By ``2-out-of-3", $\alpha ''+G'$ is a quasi-isomorphism. 

Let $\xymatrix@1{\frak k:I(\Om(\bc C^*K)^{\oplus 2})\;\;\ar @{>->} [r]&I\bigl(\Om(\bc C^*K^{\otimes 2})\bigr)}$ be the natural inclusion. The triangle  
$$\xymatrix{\Om(\bc C^*K^{\otimes 2})\cof ^{\hat\frak j_1}_{\simeq}\ar @/_2pc/ [rr]_{\frak j_1}^\simeq &(Q,d)\cof  ^(0.35){\frak j_1+ \frak k}&I\bigl(\Om(\bc C^*K^{\otimes 2})\bigr)}$$
then commutes, proving that $\frak j_1+\frak k$ is also a quasi-isomorphism.

We now have a commuting diagram
$$\xymatrix{(Q,d)\ar @{>->} [d]_\simeq^{\frak j_1+\frak k}\ar [r]^{\alpha ''+G'}_\simeq&S^*|K|^2\ar @{->>}[d]^{S^*\Delta}\\
I\bigl(\Om(\bc C^*K^{\otimes 2})\bigr)\ar [r]^(0.6){H'}&S^*|K|}$$
and therefore, since $\frak j_1+\frak k$ is a free extension and a quasi-isomorphism and $S^*\Delta$ is surjective, there is a homotopy $\xymatrix@1{\hat H:I\bigl(\Om(\bc C^*K^{\otimes 2})\bigr)\ar [r]&S^*|K|^2}$ such that  $\widehat H\circ (\frak j_1+\frak k)=\alpha ''+G'$ and $S^*\Delta \circ \widehat H=H'$.  Furthermore, $\widehat H$ is a quasi-isomorphism by ``2-out-of-3".  

Let $\hat\alpha=\widehat H\circ \frak j_2$. Then $\hat \alpha$ is a quasi-algebra quasi-isomorphism, as it is homotopic to a quasi-algebra quasi-isomorphism.  Moreover,
$$\hat\alpha\iota =\widehat H\frak k \frak j_2=G'\frak j_2=\alpha '=S^*pr_1\circ\hat\gamma+S^*pr_2\circ\hat\gamma,$$
while
$$S^*\Delta \circ \hat\alpha=H'\frak j_2=\hat \gamma\vp,$$ i.e, the following diagram commutes exactly and is therefore a noncommutative model.
$$\xymatrix{\Om(\bc C^*K\otimes \bc C^*K)\ar [d]_\simeq ^{\hat\alpha}\ar [r]^(0.65)\vp&\Om\bc C^*K\ar [d]_\simeq^{\hat\gamma}\\
S^*|K|^2\ar [r]^{S^*\Delta}&S^*|K|}$$

To complete the proof, compose $\hat\alpha$ and $\hat\gamma$ with the natural cochain algebra quasi-isomorphisms
$$\xymatrix{S^*|K|^2\ar [r]^(0.3){\simeq}&CU^*X^2\quad\text{and}\quad S^*|K|\ar [r]^(0.65){\simeq}&CU*X}$$
to obtain $\tilde\alpha$ and $\tilde\gamma$.
\qed\enddemo

The next step in our simplification of the free loop model is to find an appropriate model for the topological path fibration
$$\xymatrix{p:X^I\ar @{->>}[r]&X\times X:\ell\;\ar @{|->}[r]&(\ell(0),\ell(1))}.$$
The advantage to working directly with cubical cochains is that we are able to define an explicit model morphism, extending $\tilde\alpha$, from $\Om(\bc C^*K\otimes \bc C^*K)\odot \bc C^*K$ (defined exactly as in section 1.2) to $CU^*X^I$.  This may be possible simplicially as well, but the formulas are most probably not nearly as simple.

Given $x\in X$, let $\ell _x$ denote the constant path at $x$. Let $\xymatrix@1{i:X\ar [r]&X^I}$ be defined by $i(x)=\ell _x$. Let $\xymatrix@1{H:CU_*X^I\ar [r]&CU_{*+1}X^I}$ be the chain homotopy such that for all $T\in CU_nX^I$,
$$\xymatrix{H(T)(t_0,...,t_n):=T(t_1,...,t_n)(t_0\cdot -):I\ar [r]&X.}$$
It is easy to see that $[d,H](T)=\ell _{T(-,...,-)(0)}-T,$
i.e., $[d,H]=\iota \pi -Id$,
where $$\xymatrix@1{\iota :\im CU_*i\ar [r]&CU_*X^I}$$ is the inclusion and $$\xymatrix{\pi:CU_*X^I\ar [r]&\im CU_*i}$$ is defined by $\pi(T)=\ell _{T(-,...,-)(0)}$. Observe that $H\iota =0$, as $H(T)$ is degenerate if the path $T(t_1,...,t_n)$ is constant for all $(t_1,...,t_n)$.  The homotopy $H$ therefore induces a homotopy
$$\xymatrix{H':\coker CU_*i\ar [r]&\coker CU_{*+1}i: [T]\;\ar @{|->}[r]&[H(T)]}$$
satisfying $[d,H']=-Id,$ i.e., $H'$ is a contracting homotopy.  

Upon dualizing, we obtain a  cochain homotopy
$$\xymatrix{J=-(H')^{\sharp}:\ker CU^*i \ar [r]&\ker CU^{*-1}i}$$
such that $[d^{\sharp}, J]=Id$.  Note that $[d^{\sharp}, J]=Id$ implies that $d^{\sharp} J^2=J^2d^{\sharp}$.

It is important for the constructions in section  2 to observe as well that $J$ is a $(0,Id)$-derivation, i.e.,  
$$J(fg)=J(f)\cdot g$$
for all $f,g\in \ker CU^*i$, since $\overline\Delta H=(\iota\pi \otimes H - H\otimes Id)\overline\Delta$, which implies  that $\overline\Delta H'=-(H'\otimes Id)\overline\Delta$, where $\overline\Delta$ denotes the reduced coproduct on $CU_* X^I$ or $\coker CU_* i$. 

We resume the construction above in the following proposition.

\proclaim{Proposition 1.3.2}   Let $\xymatrix@1{i:X\ar [r]&X^I}$ be defined by $i(x)=\ell _x$, the constant path at $x$. There is a natural cochain homotopy $$\xymatrix@1{J:\ker CU^*i \ar [r]&\ker CU^{*-1}i}$$ such that 
\roster
\item $[d^\sharp, J]=Id$,
\item $d^\sharp J^2=J^2d^\sharp$, and
\item $J(fg)=J(f)\cdot g$.
\endroster\endproclaim

We now apply $J$ to giving an explicit definition of a model morphism $\tilde\beta $ from $\Om (\bc C^*K\otimes \bc C^*K)\odot \bc C^*K$ to $CU^*X^I$.

\proclaim{Proposition 1.3.3} There is a noncommutative model
$$\xymatrix{\Om(\bc C^*K\otimes \bc C^*K)\ar [d]^{\tilde\alpha}_{\simeq}\ar [r]^(0.4){\iota}&\Om (\bc C^*K\otimes \bc C^*K)\odot \bc C^*K\ar [d]_{\simeq}^{\tilde\beta}\\
CU^*(X\times X)\ar [r]^{CU^*p}&CU^*X^I}$$
where $\tilde\beta$ is defined recursively by
$$\tilde \beta (w\otimes c)=\sn {wc}\tilde \beta (1\otimes c)\cdot CU^*p\circ\tilde\alpha (w)$$
and 
$$\tilde \beta (1\otimes c)=J\tilde \beta D(1\otimes c).$$
\endproclaim

\demo{Proof}  Suppose that if $\tilde \beta$ is defined in accordance with the formulas in the statement above on $\Om (\bc C^*K\otimes \bc C^*K)\odot (\bc C^*K)_{<n}$, then $d^{\sharp}\tilde \beta =\tilde \beta D$  and $CU^*i\circ \tilde\beta =\tilde\gamma\circ\bar\vp$ when restricted to  this same complex.

Let $c\in (\bc C^*K)_n$.  It is clear from the definition of $D$ that $\tilde \gamma\bar\vp D(1\otimes c)=0$.  Thus, by the induction hypothesis, $\tilde \beta D(1\otimes c)\in \ker CU^*i$, which implies that $J\tilde \beta D(1\otimes c)$ is defined.  Moreover,
$$\align
d^{\sharp}J\tilde\beta D(1\otimes c)&=\tilde \beta D(1\otimes c)-Jd^{\sharp}\tilde\beta D(1\otimes c)\\
&=\tilde \beta D(1\otimes c)-J\tilde\beta D^2(1\otimes c)\\
&=\tilde \beta D(1\otimes c)
\endalign$$
and $CU^*iJ\tilde\beta D(1\otimes c)=0=\tilde\gamma\vp (1\otimes c)$, so we can set 
$$\tilde\beta (1\otimes c)=J\tilde\beta D(1\otimes c).$$
For the usual reasons, $\tilde\beta$ is then a quasi-algebra quasi-isomorphism.
\qed\enddemo

We can now twist together the models of Theorem 1.3.1 and of Proposition 1.3.3, obtaining a noncommutative model
$$\xymatrix{\Om\bc C^*K\ar [r]\ar [d]^{\tilde\gamma}_\simeq&\Om \bc C^*K\odot \bc C^*K\ar [d]^{\tilde\delta}_\simeq\\
CU^*X\ar [r]^{CU^*e}&CU^*\Cal LX}$$
where $\Om \bc C^*K\odot \bc C^*K$ is defined exactly as in section 1.2 and 
$$\tilde \delta (w\otimes c):=\sn {wc}\tilde\delta (1\otimes c)\cdot CU^*e\circ \tilde\gamma(w)$$
and
$$\tilde\delta (1\otimes c):=CU^*jJ\tilde\beta D(1\otimes c).$$
Recall that $\xymatrix@1{j:\Cal LX\ar [r]&X^I}$ is the natural inclusion and that $\xymatrix@1{e:\Cal LX\ar@{->>}[r]&X}$ is the basepoint evaluation. 

An easy Zeeman's comparison theorem argument shows that $\tilde\delta$ is a quasi-isomorphism.  Furthermore, applying the Eilenberg-Moore spectral sequence of algebras, one obtains that $\tilde\delta$ induces an algebra isomorphism on the $E_\infty$-terms, which is not quite as strong as saying that it is a quasi-algebra map, but is good enough for our purposes in chapter 2. By arguments similar to those in \cite {DH4}, we can show that $\tilde\delta$ is truly a quasi-algebra map in certain special cases, and it may perhaps be a quasi-algebra map in general.

\subhead 1.4 Linearization of the free loop model\endsubhead
\noindent  In this section we simplify even further the free loop model, making it as small as possible, to facilitate homotopy orbit space computations in section 2. 

Consider the surjection 
$$\xymatrix{\ve_K\otimes Id: \Om \bc C^*K\otimes \bc C^*K \ar [r]&C^*K\otimes \bc C^*K.}$$
Extend the differential on $C^*K$ to a differential $\bbd$ on $C^*K\otimes \bc C^*K$, which is a free right $C^*K$-module, by
$$\bbd (1\otimes c):=(\ve _K\otimes Id)\overline D (1\otimes c),$$
extended as a right module derivation. Consequently, $\bbd (\ve_K\otimes Id)=(\ve_K\otimes Id)\overline D$, which implies in turn that $\bbd^2=0$.

Let $\xymatrix@1{\pi :\bc C^*K\ar [r]&C^*K}$ be the projection onto (desuspended) linear terms, i.e., 
$$\pi (sx_1|\cdots |sx_n):=\left \{\aligned x_1&:n=1\\
								0&:n>1.\endaligned\right .$$
Using the notation introduced in  the definition of $\Om (\bc C^*K\otimes \bc C^*K)\odot \bc C^*K$ in section 1.2, we obtain the following explicit formula for $\bbd$, when $c=sx_1|\cdots |sx_n$.
$$\bbd(1\otimes c)=1\otimes d_\bc c + \sum _{i<j} \pi \bigl(\lambda _{ij,k}(c)\star \lambda _{ij}^k(c)\bigr)\otimes sx_i|\cdots |sx_j.$$
Notice that it is entirely possible that $\bbd (1\otimes c)\not=1\otimes d_\bc c$, since the formula for $\psi_K$ implies that if $\deg x\gg 0$, then $\psi _K(\si x)$ has a nonzero summand in $T^{>1}\si C_+K\otimes \si C_+K$.

Define now a left $C^*K$-action on $C^*K\otimes \bc C^*K$ by 
$$(x\otimes 1)\cdot (1\otimes c):=(\ve_K\otimes Id)\bigl((\si (sx) \otimes 1)\cdot (1\otimes c)\bigr)$$
for all $x\in C^*K$ and $c\in \bc C^*K$.  If $\overline D (1\otimes c)-1\otimes d_\bc c=\sum_i \si (a_i)\otimes b^i$, then 
$$\align
\bbd&\bigl((x\otimes 1)\cdot (1\otimes c)\bigr)\\
=&(\ve_K\otimes Id) \bigl((\si (sdx)\otimes 1)\cdot (1\otimes c)+\sn x (\si (sx)\otimes 1)\cdot \overline D(1\otimes c)\bigr)\\
=&(dx\otimes 1)\cdot (1\otimes c)\\
&+\sn x (\ve_K\otimes Id)\biggl((\si(sx)\otimes1)\cdot\bigl(1\otimes d_\bc c\\
&\qquad\qquad\qquad\qquad\qquad\qquad\qquad\qquad+\sum _i \sn {(a_i +1)b^i}(1\otimes b^i)(\si (sa_i)\otimes 1)\bigr)\biggr)\\
=&(dx\otimes 1)\cdot (1\otimes c)\\
&+\sn x \biggl((x\otimes1)\cdot\bigl(1\otimes d_\bc c+\sum _i \sn {(a_i +1)b^i}(1\otimes b^i)(a_i \otimes 1)\bigr)\biggr)\\
=&(dx\otimes 1)\cdot (1\otimes c)+\sn x (x\otimes1)\cdot \bbd (1\otimes c),
\endalign$$
i.e., the left $C^*K$ action commutes with the differential.

Again using the notation of section 1.2, we can write
$$(x\otimes 1)\cdot (1\otimes c)=x\otimes c - \sum _{i<j} \pi \bigl(x\star\lambda _{ij,k}(c)\star \lambda _{ij}^k(c)\bigr)\otimes sx_i|\cdots |sx_j.$$

The following proposition summarizes the observations above. Let $C^*K\tot \bc C^*K$ denote $C^*K\otimes \bc C^*K$ endowed with the differential $\bbd$ and the $C^*K$-bimodule structure defined above.

\proclaim{Proposition 1.4.1}  There is a twisted bimodule extension 
$$\xymatrix{C^*K\ar [r]&C^*K\tot \bc C^*K\ar [r]&\bc C^*K}$$
such that $\xymatrix@1{\ve _K\otimes Id :\Om\bc C^*K\tot \bc C^*K\ar [r]&C^*K\tot \bc C^*K}$ is a map of differential right $C^*K$-modules.  In particular, $\ve_K\otimes Id$ is a quasi-isomorphism.\endproclaim

Note that $\ve_K\otimes Id$ is not a bimodule map itself, since it is possible that for some $ w\in \perp ^{>1}sC^+K$ and $c\in \bc C^*K$, the product $(\si (w)\otimes 1)(1\otimes c)$ has a nonzero summand in $\si (sC^*K)\otimes \bc C^*K$, i.e., that $(\ve _K\otimes Id)\bigl((\si (w)\otimes 1)(1\otimes c)\bigr)\not=0$, even though $\si (w)\in \ker \ve _K$.  On the other hand if we filter both $\Om\bc C^*K\odot \bc C^*K$ and $C^*K\tot \bc C^*K$ by degree in the left tensor factor, then $\ve_K\otimes Id$ induces an isomorphism of bigraded bimodules on the $E_{\infty}$-terms of the associated spectral sequences, so $\ve_K\otimes Id$ is almost a quasi-bimodule map.  It may even be possible to define explicitly a cochain homotopy ensuring that $\ve_K\otimes Id$ truly is a quasi-bimodule map.

For the constructions in sections 2 and 3, we need a quasi-isomorphism $$\xymatrix@1{\Ups:C^*K\tot \bc C^*K\ar [r]^(.65)\simeq &CU^*\Cal L X},$$
which we obtain as follows. Recall that $\ve_K$ has a differential, though not multiplicative, section
$$\xymatrix{\sigma _K:C^*K\ar [r]&\Om\bc C^*K:x\;\ar @{|->}[r]&\si (sx).}$$ 
Consider the following commutative diagram of cochain complexes and maps.
$$\xymatrix{C^*K\tot \perp ^{\leq 2}sC^+K\ar[d]^(0.6){\text{incl.}}\ar [rr]^{\sigma _K\otimes Id}&&\Om\bc C^*K\odot \perp ^{\leq 2}sC^+K\ar [r]^{\text{incl.}}&\Om\bc C^*K\odot \bc C^*K\ar @{->>}[d]_\simeq^{\ve _K\otimes Id}\\
C^*K\tot \bc C^*K\ar [rrr]^{Id}&&&C^*K\tot \bc C^*K}$$
Since the inclusion map on the left is a free extension of cochain complexes and $\ve_K\otimes Id$ is a surjective quasi-isomorphism, we can extend $\sigma _K\otimes Id$ to a cochain map 
$$\xymatrix{\hat \sigma:C^*K\tot \bc C^*K\ar [r]&\Om \bc C^*K\odot \bc C^*K}$$
such that $(\ve _K\otimes Id)\hat \sigma _K=Id$, i.e., $\hat \sigma$ is a section of $\ve_K\otimes Id$. In particular, $\hat\sigma$ is a quasi-isomorphism and for all $x\otimes c\in C^*K\otimes \bc C^*K$ 
$$\si (sx)\otimes c-\hat \sigma (x\otimes c)\in \ker \ve_K\otimes \bc C^*K.\tag 1.4.1$$
Furthermore, like $\ve_K\otimes Id$, $\hat\sigma$ induces an isomorphism of bigraded bimodules on the $E_{\infty}$-terms of the usual spectral sequences and is a quasi-bimodule map if and only if $\ve_K\otimes Id$ is. 
 
We now define $\Ups$ to be the composition below.
$$\xymatrix{C^*K\tot \bc C^*K\ar [r]^(.45){\hat \sigma}\ar @/_1.5pc/ [rr]_{\Ups}&\Om \bc C^*K\odot \bc C^*K\ar [r]^(.6){\tilde \delta}&CU^*\Cal LX}$$
Observe that (1.4.1) implies that for all $x\otimes c\in C^*K\tot \bc C^*K$, 
$$\Ups (x\otimes c)=\tilde \delta (\si (sx)\otimes c)\tag 1.4.2$$
since $\ker \ve_K\subseteq \ker \tilde \gamma$.

\definition{Definition}  Let $X$ be a $1$-connected space with the homotopy type of a finite-type CW-complex, and let $K$ be a finite-type, $1$-reduced simplicial set such that $|K|\simeq X$. The twisted $C^*K$-bimodule extension 
$$fls^*(X):=C^*K\tot \bc C^*K$$
together with the quasi-isomorphism
$$\xymatrix{\Ups: fls^*(X)\ar [r]^\simeq & CU^*\Cal LX}$$
is a {\sl thin free loop model} for $X$.\enddefinition

We conclude this section with an important observation concerning the relation between $\Ups$ and the product $\psi _K^\sharp$.

\proclaim{Proposition 1.4.2} Suppose that $\psi _K^{\sharp}$, the multiplication on $\bc C^*K$, is commutative.  If $x\in C^*K$ is a cycle, then  $\Ups (1\otimes c\star sx)= \Ups (1\otimes c)\cdot \Ups (1\otimes sx)$. In particular
$$\Ups (1\otimes sx_1\star\cdots\star sx_n)=\Ups (1\otimes sx_1)\cdots\Ups(1\otimes sx_n)$$
if $dx_i=0$ for some $i$. 
\endproclaim

Consequently, if $\psi _K^{\sharp}$ is commutative, then the commutator $[\Ups (1\otimes sx),\Ups(1\otimes sy)]=0$, for all cycles $x,y\in C^*K$.

\demo{Proof} Recall from section 1.2 that if $\psi _K^{\sharp}$ is commutative, then we can choose the multiplication in the model $\Om (\bc C^*K\otimes \bc C^*K)\odot \bc C^*K$ so that $(1\otimes c)\cdot (1\otimes c')=1\otimes c\star c'$ for all $c,c'\in C^*K$. According to Proposition 1.3.3 we then have that
$$\align
\tilde \beta (1\otimes c\star c')=&J\tilde\beta D(1\otimes c\star c')\\
=&J\tilde\beta D\bigl( (1\otimes c)(1\otimes c')\bigr)\\
=&J\tilde\beta\bigl(D(1\otimes c)\cdot (1\otimes c')+\sn c (1\otimes c)\cdot D(1\otimes c')\bigr)\\
=&J\tilde\beta D(1\otimes c)\cdot \tilde\beta (1\otimes c')+\sn c J\tilde \beta (1\otimes c)\cdot \tilde\beta D(1\otimes c')\\
=&\tilde\beta (1\otimes c)\cdot \tilde\beta (1\otimes c')+\sn c J^2\tilde\beta D(1\otimes c)\cdot d^{\sharp}\tilde\beta (1\otimes c').
\endalign$$
Thus, since $d^{\sharp} J^2=J^2d^{\sharp}$,
$$\align
\Ups (1\otimes c\star c')=&CU^*j\tilde\beta (1\otimes c\star c')\\
=&\Ups (1\otimes c)\cdot \Ups (1\otimes c')+\sn c d^{\sharp}CU^*j J^2 \tilde\beta (1\otimes c)\cdot d^{\sharp}\Ups (1\otimes c)\\
=&\Ups (1\otimes c)\cdot \Ups (1\otimes c')+\sn c d^{\sharp}CU^*j J^2 \tilde\beta (1\otimes c)\cdot \Ups \bbd(1\otimes c).
\endalign$$
Thus, if $x\in C^*K$ is a cycle, then $\bbd (1\otimes sx)=0$ and so $\Ups (1\otimes c\star sx)= \Ups (1\otimes c)\cdot \Ups (1\otimes sx)$. The second part of the statement follows by induction.\qed\enddemo


\head {2. Homotopy orbit spaces} \endhead

\noindent In this section we construct a noncommutative model $hos^*(X)$ for the homotopy orbit space $(\Cal LX)_{hS^1}$ of the natural $S^1$-action on the free loop space $\Cal LX$.  The form of $hos^*(X)$ is, not surprisingly, similar to that of the complex that gives the cyclic cohomology of an algebra.  The author is grateful to Nicolas Dupont for the ideas he contributed during our discussions of  $(\Cal LX)_{hS^1}$ over the years.

We begin by proving the existence of a very special family of primitive elements in the reduced cubical chains on $S^1$ and then studying its properties.  We then introduce a particularly useful resolution of the cubical chains on $ES^1$ as a module over the cubical chains on $S^1$, which we apply  to constructing a model of the homotopy orbit space of any $S^1$-action.  In the final part of this section we specialize to the case of $\Cal LX$, obtaining a small, noncommutative model for $(\Cal LX)_{hS^1}$ as an extension of the thin free loop model $fls^*(X)$.  

\subhead 2.1 A special family of primitives\endsubhead

\noindent Let $CU_*(X)$ denote the reduced cubical chains on a topological space $X$.  We begin by defining a suspension-type degree +1 operation on $CU_*S^1$.

\definition{Definition}  Given any continuous map $\xymatrix@1{f:I^n\ar [r]&\Bbb R^+}$, let 
$$\xymatrix{\hat f:I^{n+1}\ar [r]&\Bbb R^+:(t_0,...,t_n)\ar @{|->}[r]& t_0\cdot f(t_1,...,t_n)}.$$
If $\xymatrix@1{T:I^n\ar [r]&S^1}$ is an $n$-cube such that $T(t_1,...,t_n)=e^{i2\pi f(t_1,...,t_n)}$, let $\sigma (T)$ be the $(n+1)$-cube defined by  
$$\sigma (T)(t_0,...,t_{n}):=e^{i2\pi \hat f(t_0,...,t_n)},$$
where we are considering $S^1$ as the unit circle in the complex plane, i.e.,
$$S^1=\{e^{i\theta}\mid 0\leq \theta\leq 2\pi\}.$$\enddefinition

\remark{Remark} It is clear that $\sigma (T)$ is degenerate if $T$ is degenerate.  The operation $\sigma$ can therefore be extended linearly to all of $CU_*S^1$.\endremark\medskip

As the next lemma states, $\sigma $ is a contracting homotopy in degrees greater than one and is  a $(0,Id)$-coderivation.

\proclaim{Lemma 2.1.1} Let $T\in CU_*S^1$.
\roster
\item If $\deg T\geq 2$, then $d\sigma (T)=T-\sigma (dT)$ where $d$ is the usual differential on $CU_*S^1$.
\item $\overline\Delta (\sigma (T))=\sigma (T_i)\otimes T^i$, where $\overline\Delta$ is the usual reduced coproduct on $CU_*S^1$ and $\Delta (T)=T_i\otimes T^i$ (using the Einstein summation convention).
\endroster
\endproclaim

Simple calculations, applying the definitions of the cubical differential and the cubical coproduct, as given for example in \cite {Mas} and \cite {An}, suffice to prove this lemma.

We now apply the $\sigma $ operation to the recursive construction of an important family of elements in $CU_*S^1$.

\definition{Definition} Let $\xymatrix@1{T_0:I\ar [r]&S^1}$ be the $1$-cube defined by $T_0(t)=e^{i2\pi t}$.  Given $T_k\in CU_{2k+1}S^1$ for all $k<n$, let $T_n$ be the $(2n+1)$-cubical chain defined by 
$$T_n:=\sigma \biggl(\sum _{i=1}^n T_{i-1}\cdot T_{n-i}\biggr)\in CU_{2n+1}S^1.$$
Let $\Cal T:=\{T_n\mid n\geq 0\}$.\enddefinition

\example{Examples} It is easy to see that
$$T_1(t_0, t_1, t_2)=e^{i2\pi t_0(t_1+t_2)}$$
and that 
$T_2=U+V$ where
$$U(t_0,...t_4)=e^{i2\pi t_0(t_1+(t_2+t_3)t_4)}\text{ and }V(t_0,...t_4)=e^{i2\pi t_0((t_1+t_2)t_3+t_4)}.$$
\endexample

\proclaim{Proposition 2.1.2} The family $\Cal T$ satisfies the following properties.
\roster
\item $dT_0=0$ and $0\not=[T_0]$ in $H_1S^1$.
\item $dT_n=\sum _{i=1}^n T_{i-1}T_{n-i}$ for all $n>0$.
\item Every $T_n$ is primitive in $CU_*S^1$
\endroster
\endproclaim

\demo{Proof} Points (1) and (2) are easy consequences of Lemma 2.1.1.  It is well known that $T_0$ represents the unique nonzero homology generator of $H_*S^1$.

An easy inductive argument applying Lemma 2.1.1(2) proves point (3), since if $T_k$ is primitive for all $k<n$, then the sum $\sum _{i=1}^n T_{i-1}\cdot T_{n-i}$ is also primitive, as it is symmetric and all factors are of odd degree.\qed\enddemo

Let $<\Cal T>$ denote the subalgebra of $CU_*S^1$ generated by the family $\Cal T$.  Since all the $T_n$'s are primitive, $<\Cal T>$ is a sub Hopf algebra of $CU_*S^1$. Proposition 2.1.2(1) and (2) imply that $<\Cal T>$ is closed under the differential, and that the inclusion 
$$\xymatrix{<\Cal T>\ar @{^{(}->}[r]^\simeq&CU_*S^1}$$
is a quasi-isomorphism.

\subhead 2.2 A useful resolution of $CU_*ES^1$\endsubhead

\noindent We now put the family $\Cal T$ to work in constructing a simple, neat resolution of $CU_*ES^1$ as a $CU_*S^1$-module. To understand why this is important, recall that a special case of Moore's theorem (cf., \cite {Mc}, Thm. 7.27) states that for any left $S^1$-space $X$ and any free $CU^*S^1$-resolution $(Q,d)$ of $CU_*ES^1$, there is a diagram of quasi-isomorphisms of chain complexes
$$\xymatrix{(Q,d)\otimes _{CU_*S^1} CU_*X\ar [r]^(0.6)\simeq\ar [d] ^{\pi} &CU_*(X_{hS^1})\ar[d]^\pi \\
			(Q,d)\otimes _{CU_*S^1} \Bbb Z \ar [r]^(0.6) \simeq &CU_*(BS^1)}$$
where the projection maps $\pi$ are induced by the map $\xymatrix@1{X\ar [r]& {*}}$.
Hence, 
$$\tor _*^{CU_*S^1}(CU_*ES^1, CU_*X)\cong H_*(X_{hS^1}),$$
and so a resolution of $CU_*ES^1$ provides us with a general tool for computing $H_*(X_{hS^1})$ for an arbitrary $S^1$-space $X$. 

Let $\Gamma $ denote the divided powers algebra functor.  Recall that if $w$ is in even degree, then
$$\Gamma w=\bigoplus _{k\geq 0}\Bbb Z\cdot w(k),$$
where $\deg w(k)=k\cdot |w|$, $w(0)=1$, $w(1)=w$ and $w(k)w(l)=\binom {k+l}{k}w(k+l)$.  Furthermore, $\Gamma w$ is in fact a Hopf algebra, where the coproduct is specified by $\Delta (w)=w\otimes 1 +1\otimes w$. 

Consider $(\Lambda u,0)=(\Bbb Z\cdot u,0)$, where $u$ is in degree $1$, and its acyclic extension $(\Gamma v\otimes \Lambda u, \del )$, where $v$ is in degree $2$ and $\del v(k):=v(k-1)\otimes u$ for all $k\geq 1$.  There is a chain algebra quasi-isomorphism
$$\xymatrix{\xi:CU_*S^1 \ar [r]^{\simeq}&(\Lambda u,0)}$$
defined by $\xi (T_0)=u$ and $\xi (T)=0$ for all other cubes $T$.

Define a semifree extension of right $CU_*S^1$-modules 
$$\xymatrix{\iota :CU_*S^1\ar [r]&\bigl(\Gamma v\otimes CU_*S^1, \tilde\del)}$$
by 
$$\tilde \del (v(n)\otimes 1):=\sum _{k=0}^{n-1} v(n-k-1)\otimes T_k.$$
It is an immediate consequence of Proposition 2.1.2(2) that $\tilde\del ^2=0$.  Furthermore, for all $n\geq 1$,
$$(Id\otimes \xi)\tilde \del (v(n)\otimes 1)=v(n-1)\otimes u=\del (Id\otimes \xi )(v(n)\otimes 1),$$
which implies that the $CU_*S^1$-module map
$$\xymatrix{Id\otimes \xi: (\Gamma v\otimes CU_*S^1 ,\tilde\del)\ar [r]&(\Gamma v\otimes \Lambda u,\del)}$$
is a differential map.  A quick Zeeman's Comparison Theorem argument then shows that $Id\otimes\xi$ is a quasi-isomorphism, so that $ (\Gamma v \otimes CU_*S^1, \tilde\del)$ is acyclic.

We claim that $(\Gamma v \otimes CU_*S^1, \tilde\del)$ is a $CU_*S^1$-resolution of $CU_*ES^1$.  To verify this, observe that there is a commutative diagram of $CU_*S^1$-modules
$$\xymatrix{CU_*S^1\ar [d]^{\iota}\ar [r]^{CU_*j}& CU_*ES^1\ar @{->>}[d]_\simeq\\
(\Gamma v \otimes CU_*S^1, \tilde\del)\ar [r]^(.7){\simeq}&\Bbb Z}$$
where $j$ is the inclusion of $S^1$ as the base of the construction of $ES^1$, which is an $S^1$-equivariant map.  Since $(\Gamma v \otimes CU_*S^1, \tilde\del)$ is a semifree extension and the map from $CU_*ES^1$ to $\Bbb Z$ is a surjective quasi-isomorphism, we can extend $CU_*j$ to a $CU_*S^1$-module map
$$\xymatrix{\varepsilon:(\Gamma v \otimes CU_*S^1, \tilde\del)\ar [r]&CU_*ES^1}$$
which is a quasi-isomorphism by ``2-out-of-3".

\subhead 2.3 Modeling $S^1$-homotopy orbits\endsubhead

\noindent Let $X$ be any (left) $S^1$-space, where $\xymatrix@1{g:S^1\times X\ar [r]&X}$ is the action map.  There is then a natural $CU_*S^1$ module structure on $CU_*X$, given by the composition
$$\xymatrix{CU_*S^1\otimes CU_*X\ar [r]^{EZ}_\simeq\ar@/_1.5pc/ [rr]_\kappa&CU_*(S^1\times X)\ar [r]^(0.6){CU_*g}&CU_*X}$$
where $EZ$ denotes the Eilenberg-Zilber equivalence. Observe that $\kappa$ is a coalgebra map, as it is the composition of two coalgebra maps.  Since Moore's Theorem implies that 
$$H_*(X_{hS^1})\cong H_*\bigl((\Gamma v \otimes CU_*S^1, \tilde\del)\underset CU_*S^1\to\otimes CU_*X \bigr),$$
we need to try to understand better the complex $(\Gamma v \otimes CU_*S^1, \tilde\del)\underset CU_*S^1\to\otimes CU_*X $.

Define an extension $( \Gamma v\otimes CU_*X, D)$ of $\Gamma v$ by
$$D(v(n)\otimes U):=v(n)\otimes dU+\sum _{k=1}^{n-1} v(n-k-1)\otimes\kappa (T_k\otimes U).$$
We again use Proposition 2.1.2(2) to verify that $D^2=0$. Observe that $(\Gamma v\otimes CU_*X, D)$ is naturally a chain coalgebra that is a cofree left $\Gamma v$-comodule, since $\kappa $ is a coalgebra map and each $T_k$ is primitive.

It is then easy to show that the following two maps are chain isomorphisms, one inverse to the other.
$$\xymatrix{(\Gamma v \otimes CU_*S^1, \tilde\del)\underset CU_*S^1\to\otimes CU_*X\ar [r]&( \Gamma v\otimes CU_*X, D)\\
v(n)\otimes U\otimes V\;\ar @{|->}[r]&v(n)\otimes \kappa (U\otimes V)}$$

$$\xymatrix{( \Gamma v\otimes CU_*X, D)\ar [r]&(\Gamma v \otimes CU_*S^1, \tilde\del)\underset CU_*S^1\to\otimes CU_*X\\
v(n)\otimes V\;'\ar @{|->}[r]&v(n)\otimes 1\otimes V'}$$
Thus, $H_*(X_{hS^1})\cong H_*( \Gamma v\otimes CU_*X, D)$.

In this chapter we are interested in cohomology calculations and so must dualize this model.  Dualizing $\kappa$ directly poses a problem, however, since
$$(CU_*S^1\otimes CU_*X)^{\sharp}\not\cong CU^*S^1\otimes CU^*X$$
because the cubical chain complex on a space is not of finite type.  We can avoid this problem by observing that it is enough to dualize the composition
$$\xymatrix{<\Cal T>\otimes CU_*X\ar @{^{(}->}[r]^\iota_{\simeq}&CU_*S^1\otimes CU_*X\ar [r]^(.6)\kappa&CU_*X.}$$
Let $\xymatrix@1{j_n: \Bbb Z\cdot T_n\otimes CU_*X\ar [r]& <\Cal T>\otimes CU_*X}$ denote the natural inclusion of graded modules.  Let $T_n^\sharp\in CU^{2n+1}S^1$ denote the cochain such that $T_n^\sharp (T_n)=1$ and $T_n^\sharp (T)=0$ if $T$ is any other $(2n+1)$-cube. Let $<\Cal T>^\sharp =\operatorname{Hom}(<\Cal T>, \Bbb Z)$

For each $n\geq 0$, define a linear map $\xymatrix@1{\omega _n: CU^*X\ar[r]&CU^{*-(2n+1)}X}$ of degree $-(2n+1)$ by
$$j_n^\sharp\circ (\kappa\iota)^{\sharp}\bigl(f):=T_n^\sharp\otimes \omega _n(f),$$
where
$$\xymatrix{CU^*X\ar [r]^(0.4){(\kappa\iota)^\sharp}&<\Cal T>^\sharp\otimes CU^*X\ar [r]^(.475){j_n^\sharp}&\Bbb Z\cdot T_n^\sharp\otimes CU^*X.}$$ In other words, $\omega _n$ is the dual of $\kappa (T_n\otimes -)$.

Let $(\Lambda \upsilon\otimes CU^*X, D^\sharp)$ denote the $\Bbb Z$-dual of $( \Gamma v\otimes CU_*X, D)$. In particular $\upsilon(v)=1$. Since it is the dual of a cofree comodule,  $(\Lambda \upsilon\otimes CU^*X, D^\sharp)$ is a free, right $\Lambda \upsilon$-module. We need to identify $D^\sharp$ as precisely as possible, since
$$\H^*(\Lambda \upsilon\otimes CU^*X, D^\sharp)\cong \H^*(X_{hS^1}).$$
A simple dualization calculation gives us the following result.

\proclaim{Lemma 2.3.1} If $f\in CU^mX$, then 
$$D^\sharp(\upsilon^n\otimes f)=\upsilon ^n \otimes d^\sharp f+\sum _{k=0}^{\lceil \frac {m-1}2\rceil } \upsilon ^{n+k+1}\otimes\omega _k(f)$$
where $d^\sharp$ denotes the differential of $CU^*X$.
\endproclaim 

As a consequence of this description of $D^{\sharp}$ we obtain the following useful properties of the operators $\omega _k$.

\proclaim{Corollary 2.3.2} The operators $\omega _n$ satisfy the following properties.
\roster
\item  For all $n\geq 1$, $[d^\sharp, \omega _n]=-\sum _{k=0}^{n-1}\omega _k\circ \omega _{n-k-1}$, while $d^\sharp \omega _0=-\omega _0d^\sharp$.
\item Each $\omega _n$ is a derivation, i.e.,  $\omega _n(f\cdot g)= \omega _n(f)\cdot g +\sn {f}f\cdot \omega _n(g)$.
\endroster\endproclaim

\demo{Proof}  The proof of (1) proceeds by expansion of the equation $0=(D^\sharp)^2(1\otimes f)$.  To prove (2), expand the equation 
$$D^{\sharp} (1\otimes f\cdot g)=D^{\sharp} (1\otimes f)\cdot (1\otimes g) +(-1)^f(1\otimes f)\cdot D^{\sharp}(1\otimes g).$$
The differential $D^{\sharp}$ is a derivation, since it is the dual of the differential of a chain coalgebra.
\qed\enddemo 

\remark{Remark} This corollary implies that $\omega _0$ induces a derivation of degree $-1$ 
$$\xymatrix@1{\varpi: \H ^*X\ar [r]&\H^{*-1}X}$$
such that $\varpi^2=0$
\endremark

\subhead 2.4 The case of the free loop space\endsubhead

\noindent Let $K$ be a finite-type, $1$-reduced simplicial set such that $|K|$ has the same homotopy type as $X$.  As we saw in section 1.4, there is commutative diagram
$$\xymatrix{C^*K\ar [d]^{\widetilde\gamma \sigma _K}\ar [r]^(0.4){\iota}&C^*K\tot \bc C^*K\ar [d]_\simeq^{\Upsilon} \\
CU^*X\ar [r]^{CU^*e}&CU^*\Cal LX}$$
in which $\widetilde\gamma \sigma _K$ is a quasi-algebra quasi-isomorphism, $\iota $ is a twisted bimodule extension and $\Ups$ is a quasi-isomorphism inducing an isomorphism on the $E_{\infty}$-terms of the Eilenberg-Moore spectral sequence.  

Our goal here is to combine this thin free loop space model with the general homotopy orbit space model of the section 2.3, obtaining an extension $(\Lambda \upsilon\otimes C^*K\tot\bc C^*K,\tD)$ of $(\Lambda \upsilon, 0)$ by $fls^*(X)$, together with a quasi-isomorphism
$$\xymatrix{{\wU}:(\Lambda \upsilon\otimes C^*K\tot\bc C^*K,\tD)\ar [r]&(\Lambda \upsilon\otimes CU^*\Cal LX, D^\sharp)}$$
such that
$$\xymatrix{(\Lambda \upsilon,0)\ar @{=} [d]\ar [r]^(0.35){incl.}&(\Lambda \upsilon\otimes C^*K\tot \bc C^*K, \tD)\ar [r]^(0.6)\pi\ar [d]^{\wU}& C^*K\tot\bc C^*K\ar [d]^\Ups\\
(\Lambda \upsilon,0)\ar [r]^(0.35){incl.}&(\Lambda \upsilon\otimes CU^*\Cal LX,D^\sharp)\ar [r] ^(0.6)\pi &CU^*\Cal L X}\tag 2.4.1$$
commutes, where $\pi$ denotes the obvious projections.

We begin by an easy, though crucial, observation concerning the relations between $\Ups$ and the operations $\omega _k$.

\proclaim {Lemma  2.4.1}  For all cocycles $x\in C^*K$, $\Ups (1\otimes sx)=\omega _0\Ups (x\otimes 1)$.\endproclaim

\demo{Proof} From the definition of $\tilde\beta$ from Proposition 1.3.3 and of $D$ from section 1.2, we can show that
$$\align
\Ups (1\otimes sx)=&CU^*j\tilde\beta (1\otimes sx)\\
=&CU^*jJ\tilde \beta D(1\otimes sx)\\
=&CU^*j J\tilde \beta (\si (sx\otimes 1)-\si (1\otimes sx))\\
=&CU^*j JCU^*p\bigl( C^*pr_1 -C^*pr_2\bigr)\tilde\gamma (\si(sx)).
\endalign$$
A straightforward computation suffices to establish that $$CU^*j JCU^*p\bigl( C^*pr_1 -C^*pr_2\bigr)=\omega _0CU^*e,$$ implying that
$$\Ups (1\otimes sx)=\omega _0 CU^*e \tilde\gamma (\si (sx))=\omega _0\Ups (x\otimes 1).\qed$$
\enddemo

\definition{Restriction} Henceforth, to simplify the presentation, we assume that $\psi _K^{\sharp}$, the multiplication on $\bc C^*K$, is such that the primitives of $\bc C^*K$, i.e., the elements of $sC^+K$, are all indecomposable. \enddefinition

This is certainly a strong hypothesis, but it still allows us to treat a number of interesting cases, such as wedges of spheres.  More general cases are treated in \cite {H2}. 

The special properties of the free loop space model in the restricted case we consider are summarized in the following lemma.

\proclaim{Lemma 2.4.2} If the primitives of $\bc C^*K$ are all indecomposable, then the following properties hold.
\roster
\item The multiplication on $\bc C^*K$ is the shuffle product, which is commutative.
\item The graded algebra $C^*K$ is commutative. 
\item For all $y, x_1,...,x_n\in C^*K$,
$$\align
\bbd (y\otimes sx_1|\cdots|sx_n)= &dy\otimes sx_1|\cdots|sx_n +\sn y y\otimes d_\bc(sx_1|\cdots|sx_n)\\
&+ \sn y\bigl[yx_1\otimes sx_2|\cdots|sx_n  -\sn N yx_n\otimes sx_1|\cdots |sx_{n-1}\bigr], \endalign$$
where $N=(1+\deg x_n )(n-1 +\sum _{1\leq j<n}\deg x_j)$, and
$$(y\otimes 1)(1\otimes sx_1|\cdots|sx_n)=y\otimes sx_1|\cdots |sx_n.$$\endroster\endproclaim

In other words, the differential of $C^*K\tot \bc C^*K$ is exactly that of the usual Hochschild complex on $C^*K$, while the left $C^*K$-action is untwisted, when the primitives of $\bc C^*K$ are all indecomposable.

\demo{Proof} (1)  This is obvious.
\medskip
\noindent (2) Recall from section 1.2 that if $\psi_K ^{\sharp}$ is commutative, then $$(1\otimes c)(1\otimes c')=1\otimes c\star c'.$$ Thus,  if $x\in C^{m+1}K, y\in C^{n+1}K$, then
$$\align
\bbd (1\otimes sx\star sy)=&\bbd (1\otimes sx)\cdot (1\otimes sy)+\sn m (1\otimes sx)\cdot \bbd (1\otimes sy)\\
=&-(1\otimes s(dx))\cdot (1\otimes sy)-\sn m (1\otimes sx)\cdot(1\otimes s(dy)),
\endalign$$
whenever $\psi_K ^{\sharp}$ is commutative.
If, moreover, all primitives of $\bc C^*K$ are indecomposable, then 
$$\align
\bbd (1\otimes sx\star sy)=&\bbd (1\otimes sx|sy+\sn {mn}1\otimes sy|sx)\\
=&x\otimes sy-\sn {mn}y\otimes sx+\sn m 1\otimes s(xy)\\
&-1\otimes s(dx)|sy-\sn m1\otimes sx|s(dy)\\
&+\sn {mn}\bigl[y\otimes sx-\sn {mn}x\otimes sy+\sn n 1\otimes s(yx)\\
&\qquad\qquad\quad-1\otimes s(dy)|sx-\sn n1\otimes sy|s(dx)\bigr]\\
=&\sn m 1\otimes s\bigl([x,y]\bigr )-1\otimes s(dx)\star sy-\sn m1\otimes sx\star s(dy)\endalign$$
and so $[x,y]=0$.  Hence, $C^*K$ is commutative if all primitives of $\bc C^*K$ are indecomposable.
\medskip

\noindent (3) When all primitives of $\bc C^*K$ are indecomposable, the formulas of section 1.4 obviously reduce to those given in the statement.\qed\enddemo

We now define the desired extension 
$$\xymatrix{(\Lambda \upsilon ,0)\ar [r]&(\Lambda \upsilon\otimes C^*K\otimes\bc C^*K,\tD)\ar[r]& C^*K\odot \bc C^*K}$$ 
and show that 
$$\H^*(\Lambda \upsilon\otimes C^*K\otimes \bc C^*K, \tD)\cong \H^*(X_{hS^1}).$$
We define the extension  by
$$\tD= Id\otimes \bbd + \upsilon \cdot -\otimes S$$
where 
\roster
\item $S(y\otimes 1)=1\otimes sy$ for all $y\in C^+K$;
\item $S(1\otimes c)=0$ for all $c\in \bc C^*K$; and
\item $S(y\otimes sx_1|\cdots|sx_n)=\sum _{ j=1}^{n+1}\pm 1\otimes sx_j|\cdots|sx_n| sy|sx_1|\cdot |sx_{j-1},$ where the sign is chosen in accord with the Koszul convention (cf., section 0.1).
\endroster
Thus, for example, if $x,y,z\in C^*K$ are of degrees $l+1$, $m+1$ and $n+1$, respectively, then $S(x\otimes sy)=1\otimes (sx|sy+\sn {lm} sy|sx)=1\otimes sx\star sy$ and
$$S(x\otimes sy|sz)=1\otimes (sx|sy|sz+\sn{(m+l)n} sz|sx|sy +\sn {l(m+n)}sy|sz|sx).$$

It is obvious that $S^2=0$.  A tedious, though not difficult, combinatorial calculation, shows that $\bbd S=-S\bbd$.  The proof of this equality depends strongly on the fact that $C^*K$ is commutative; in the general case we need to add terms to $S(x\otimes c) $ to kill certain commutators \cite {H2}. Thus $\tD^2=0$, i.e., $(\Lambda \upsilon\otimes C^*K\odot \bc C^*K, \tD)$ is a cochain complex.  Indeed this is exactly the negative cyclic complex of the commutative algebra $C^*K$, looked at as a cochain complex in positive degrees, rather than as a chain complex in negative degrees.

As Jones proved in \cite {J}, $\H^*(\Cal LX_{hS^1})$ is isomorphic to the negative cyclic homology of the algebra $S^*X$, and therefore to that of $C^*K$, if $|K|\simeq X$.  Thus
$$\H^*(\Lambda \upsilon \otimes C^*K\tot \bc C^*K, \tD)\cong \H^*(\Lambda \upsilon\otimes C^*\Cal L X,D^{\sharp})\cong \H^*(\Cal LX_{hS^1}),$$
as desired. To build our model for topological cyclic homology, however,  we need a cochain quasi-isomorphism $\wU$ lifting $\Ups$ and inducing this cohomology isomorphism.  In the next theorem, we prove the existence of $\wU$ when $K$ is an odd-dimensional sphere.  Using results from \cite P, we can generalize this theorem to wedges of odd spheres and, when working over $\Bbb F_2$, to wedges of even spheres.  The essential ideas of the general proof are already present in the proof for a single odd sphere, so we restrict to this case, to simplify the presentation. In \cite {H2}, we prove the existence of $\wU$ for a somewhat larger class of spaces. 

Before stating and proving the theorem, we analyze carefully the $S^1$- action on $\Cal L S^{2n+1}$.  It is well known that $\H^*(\Cal LS^{2n+1})$ is isomorphic to the tensor product of an exterior algebra $\Lambda x$ on a generator of degree $2n+1$  with the divided powers algebra $\Gamma y$ on a generator of degree $2n$ (cf., e.g., \cite {Sm}).  The generator $x$ is represented by $CU^*e(\zeta)$, where $\zeta\in CU^{2n+1}S^{2n+1}$ represents the fundamental class of $S^{2n+1}$.  More explicitly, if $\xymatrix{U:I^{2n+1}\ar [r]&S^{2n+1}}$ is a $(2n+1)$-cube collapsing $\del I^{2n+1}$ to a  point, then
$$\xymatrix{\zeta :CU_{2n+1}S^{2n+1}\ar [r]&\Bbb Z}$$
is specified by $\zeta (U)=1$ and $\zeta (V)=0$ if $V$ is any other $(2n+1)$-cube.  Consider the transpose of $U$
$$\xymatrix{U^\flat :I^{2n}\ar [r]& (S^{2n+1})^I: (t_1,...,t_{2n})\;\ar @{|->}[r]&U(-,t_1,...,t_{2n})}.$$
Since $U$ collapses the boundary of the cube, $U^\flat(t_1,..,t_{2n})$ is always a (based) loop.  Let $\xi\in CU^{2n}\Cal LX$ be the cochain such that  $\xi (U^\flat)=1$ and $\xi (V)=0$ for any other $2n$-cube $V$.  The generator $y$ is represented by $\xi$.

Recall the definition of  $\varpi$ at the end of section 2.3. A simple, explicit calculation shows that $\varpi(x)=y,$ which implies that $\varpi(y)=0$, whence  $$\varpi(x\otimes y(m))=(m+1)\cdot y(m+1) \quad \text{and}\quad \varpi(1\otimes y(m))=0$$ 
for all $m\geq 0$, since $\varpi$ is a derivation.  In particular, $\varpi(\H^{\text{even}}\Cal LS^{2n+1})=0$, and $\xymatrix{\varpi:\H^{\text{odd}}\Cal L S^{2n+1}\ar [r]&\H^{\text{even}}\Cal L S^{2n+1}}$ is an isomorphism.

\proclaim{Theorem 2.4.3}  Let $K$ be the simplicial model of $S^{2n+1}$ with exactly two nondegenerate simplices, in degrees $0$ and $2n+1$, where $n>0$. There is a quasi-isomorphism $\xymatrix@1{{\wU}: (\Lambda \upsilon\otimes C^*K\tot \bc C^*K, \tD)\ar [r]&(\Lambda \upsilon\otimes CU^*\Cal LX,D^\sharp)}$ such that 
$$\xymatrix{(\Lambda \upsilon,0)\ar @{=} [d]\ar [r]^(0.35){incl.}&(\Lambda \upsilon\otimes C^*K\tot \bc C^*K, \tD)\ar [r]^(0.6)\pi\ar [d]^{\wU}& C^*K\tot\bc C^*K\ar [d]^\Ups\\
(\Lambda \upsilon,0)\ar [r]^(0.35){incl.}&(\Lambda \upsilon\otimes CU^*\Cal LS^{2n+1},D^\sharp)\ar [r] ^(0.6)\pi &CU^*\Cal L S^{2n+1}}$$
commutes, where $\pi $ denotes the obvious projection maps.\endproclaim

\demo{Proof} In this case, $C^*K=\Lambda z$, an exterior algebra on an odd generator of degree $2n+1$ and $\bc C^*K=\Gamma (sz)$, the divided powers algebra on an even generator of degree $2n$.  Furthermore, $\bbd=0$, i.e., $(C^*K\tot \bc C^*K,\bbd)=(\Lambda z\otimes\Gamma sz,0)$. 

We need to define a cochain map $\wU=\sum _{k\geq 0} \upsilon ^k\otimes \Ups _k$, where $\Ups_0=\Ups$ and, for all $k$, $\xymatrix@1{\Ups _k:C^*K\tot \bc C^*K\ar [r]& CU^{*-2k}\Cal LS^{2n+1}}$ is a linear map of degree $-2k$.  Since $\tD =1\otimes \bbd + \upsilon\otimes S$ and $D^\sharp = 1\otimes d^\sharp + \sum _{k\geq 0} \upsilon ^{k+1}\otimes \omega _k$, the equation $D^\sharp \wU=\wU \tD$ is equivalent to the set of equations
$$ d^\sharp \Ups_k +\sum _{i+j=k-1} \omega _i\Ups _j= \Ups _k\bbd +\Ups _{k-1}S
\tag "$(2.4.2)_k$"$$
for $k\geq 0$.
Thus, to define $\wU$, we can build up the family of $\Ups_k$'s by induction on both $k$ and wordlength in $\bc C^*K=\Gamma sx$.

For $k=0$, the equation above becomes simply $d^\sharp \Ups _0=\Ups _0\bbd $, which holds because $\Ups$ is a cochain map.

When $k=1$, the appropriate equation is 
$$d^\sharp \Ups _1=\Ups _1\bbd +\Ups _0S-\omega _0\Ups _0\tag "$(2.4.2)_1$"$$
Applied to $z\otimes 1\in C^*K\tot \bc C^*K$, the right-hand side of this equation becomes
$$0+\Ups (1\otimes sz)-\omega _0\Ups (z\otimes 1),$$
which is $0$, by Lemma 2.4.1.   Thus, we can set 
$$\Ups _1(z\otimes 1)=0.$$

When we apply the right-hand side of $(2.4.2)_1$ to $1\otimes sz$, we obtain
$$0+0-\omega_0\Ups (1\otimes sz),$$
which is equal to $-\omega _0^2\Ups (z\otimes 1)$, by Lemma 2.4.1.
We can therefore choose 
$$\Ups _1(1\otimes sz)=\omega _1\Ups (z\otimes 1),$$
by Corollary 2.3.2 (1).  By a similar argument, we can choose 
$$\Ups _1(z\otimes sz)=\omega _1\Ups (z\otimes 1)\cdot \Ups (z\otimes 1).$$ 

Suppose that $\Ups_1$ has been defined on $\Lambda z\otimes \Gamma ^{\leq l-1}(sz)$ satisfying equation $(2.4.2)_1$, where $l\geq 2$.  Applying the right-hand side of $(2.4.2)_1$ to $1\otimes sz(l)$, we obtain
$$0+0-\omega_0\Ups (1\otimes sz(l)),$$
which is a cycle of odd degree.  Since $\varpi$ is injective on odd cohomology, either $\omega _0(1\otimes sz(l))$ is a boundary or  $\omega _0^2(1\otimes sz(l))$ is not a boundary.   The second option is impossible, since $\omega _0^2(1\otimes sz(l))=-d^\sharp \omega _1(1\otimes sz(l))$, and so $\omega _0(1\otimes sz(l))$ must be a boundary, i.e., there is possible choice of $\Ups_1(1\otimes sz(l))$ satisfying $(2.4.2)_1$.  If we then set 
$$\Ups _1(z\otimes sz(l))=\Ups_1 (1\otimes sz(l))\cdot \Ups (z\otimes 1),$$ 
then equation $(2.4.2)_1$  is satisfied on $\Lambda z\otimes \Gamma ^{\leq l}(sz)$.

Suppose now that for all $k<m$, there is a linear map $\xymatrix@1{\Ups_k:\Lambda z\otimes \Gamma sz\ar [r] &CU^{*-2k}\Cal LS^{2n+1}}$ satisfying $(2.4.2)_k$ and such that 
$$\Ups _k(z\otimes 1)=0\quad\text{and}\quad\Ups _k(1\otimes sz)= \omega _k\Ups (z\otimes 1).$$  

For $k=m$, the equation we must satisfy is
$$ d^\sharp \Ups_m =\Ups _m\bbd +\Ups _{m-1}S-\sum _{i+j=m-1} \omega _i\Ups _j. 
\tag "$(2.4.2)_k$"$$
Applied to $z\otimes 1$, the right-hand side of the equation becomes
$$0+\Ups _{m-1}(1\otimes sz)-\sum _{i+j=m-1} \omega _i\Ups _j (z\otimes 1),$$
which is zero, by the induction hypotheses.   We can therefore set $\Ups _m(z\otimes 1)=0$.

If the right-hand side of the equation $(2.4.2)_m$ is evaluated on $1\otimes sz$, it becomes
$$\align
0+0-\sum _{i+j=m-1} \omega _i\Ups _j(1\otimes sz)=&-\sum _{i+j=m-1} \omega _i\omega _j\Ups (z\otimes 1)\\
=&d^\sharp \omega _m \Ups (z\otimes 1),\endalign$$
implying that we may set $\Ups (1\otimes sz)=\omega _m\Ups (z\otimes 1)$.

Suppose that $\Ups_m$ has been defined on $\Lambda z\otimes \Gamma ^{\leq l-1}(sz)$ satisfying equation $(2.4.2)_m$, where $l\geq 2$.  Applying the right-hand side of $(2.4.2)_m$ to $1\otimes sz(l)$, we obtain
$$0+0-\sum _{i+j=m-1} \omega _i\Ups _(1\otimes sz(l)),$$
which is a cycle of odd degree.  Since $\varpi$ is injective on odd cohomology, either $\sum _{i+j=m-1} \omega _i\Ups _j(1\otimes sx(l))$ is a boundary or  its image under $\omega _0$ is not a boundary.   Observe however that
 $$\align
 d^\sharp \biggl(\sum _{i+j=m-1}& \omega _{i+1}\Ups _j(1\otimes sz(l))\biggr)\\=&\sum \Sb{s+t=i}\\ i+j=m-1\endSb -\omega_s\omega _t\Ups _j(1\otimes sz(l))
 +\sum \Sb p+q=j-1\\ i+j=m-1\endSb \omega _{i+1}\omega _p\Ups _q(1\otimes sz(l)\\
 =&\sum _{i+j+k=m-1}-\omega _i\omega _j\Ups _k(1\otimes sz(l))
 +\sum \Sb i+j+k=m-1\\ i\geq 1\endSb\omega _i\omega _j\Ups _k(1\otimes sz(l))\\
 =&-\sum _{i+j=m-1} \omega _0\omega _i\Ups _j(1\otimes sz(l))\\
 =& -\omega _0\biggl(\sum _{i+j=m-1} \omega _i\Ups _j(1\otimes sz(l))\biggr),\endalign$$
and so $\sum _{i+j=m-1} \omega _i\Ups _j(1\otimes sz(l))$ must be a boundary. Hence, there is possible choice of $\Ups_m(1\otimes sz(l))$ satisfying $(2.4.2)_m$.  If we then set 
$$\Ups _m(x\otimes sz(l))=\Ups_m (1\otimes sz(l))\cdot \Ups (z\otimes 1),$$ 
then equation $(2.4.2)_m$  is satisfied on $\Lambda z\otimes \Gamma ^{\leq l}(sz)$.
\qed \enddemo

\definition{Definition}  Let $X$ be a $1$-connected space with the homotopy type of a finite-type CW-complex, and let $\xymatrix@1{\Ups: C^*K\tot \bc C^*K\ar [r]^(0.6)\simeq & CU^*\Cal LX}$ be a thin free loop model for $X$ such that all primitives of $\bc C^*K$ are indecomposable.  The twisted bimodule extension
$$hos ^*(X)=(\Lambda \upsilon\otimes C^*K\tot \bc C^*K, \tD)$$
together with the quasi-isomorphism 
$$\xymatrix{{\wU}:hos ^*(X))\ar [r]&(\Lambda \upsilon\otimes CU^*\Cal LX,D^\sharp)}$$ 
such that diagram (2.4.1) commutes, when it exists, is a {\sl thin model} of $\Cal L X_{hS^1}$.\enddefinition

\remark{Related work}  B\"okstedt and Ottosen have recently developed an approach to Borel cohomologly calculations for free loop spaces that is Eckmann-Hilton dual to the approach considered here and thus complementary to our methods \cite {BO2}.  They have constructed a Bousfield-type spectral sequence that converges to the cohomology of $(\Cal LX)_{hS^1}$.  While our model is easiest to deal with for spaces with few cells, the elementary cases for their model are Eilenberg-MacLane spaces.
\endremark


\head 3. A model for mod 2 topological cyclic homology\endhead

\noindent We begin this section by supplying the final piece of the machine with which we build a model of $TC(X;2)$: a model of the \pth-power map, for $p=2$.  We then use the machine to obtain an explicit and precise description of $tc^*(X)$.  To conclude we illustrate the power of both the $tc^*(X)$ and the $hos^*(X)$ models, by applying them to computing $\H^*(\Cal LS^{2n+1}_{hS^1})$ and $\H^*(TC(S^{2n+1};2);\Bbb F_2)$.

\subhead 3.1 The \pth- power map\endsubhead
 
\noindent The \pth-power map, $\lambda ^p$, on a free loop space $\Cal LX$ sends any loop to the loop that covers the same image $p$ times, turning $p$ times as fast, i.e., for all $\ell\in \Cal LX$ and for all $z\in S^1$
 $$\lambda^p(\ell)(z):=\ell (z^p),$$
 where we see $S^1$ as the set of complex numbers of norm $1$.
 
 There is another useful way to define $\lambda ^p$.  Let $\Cal LX^{(p)}$ denote the pullback of the iterated diagonal $\xymatrix@1{\Delta ^{(p-1)}: X\ar [r]&X^p}$ and of $\xymatrix@1{e^p:(\Cal LX)^p\ar [r]&X^p }$, i.e., the elements of $\Cal LX^{(p)}$ are sequences of loops $(\ell _1, ...,\ell _p)$ such that $\ell _i(1)=\ell _j(1)$ for all $i,j$.  Let $\xymatrix{e^{(p)}:  \Cal LX^{(p)}\ar [r]&X}$ denote the map sending a sequence of loops to their common basepoint.
 
The iterated diagonal map on $\Cal LX$ corestricts to $\xymatrix@1{\Delta ^{(p-1)}:\Cal LX\ar [r]&\Cal LX^{(p)}}$, while concatenation of loops defines a map $\xymatrix@1{\mu ^{(p-1)}: \Cal LX^{(p)}\ar [r]&\Cal LX}$, restricting to the usual iterated multiplication on $\Om X$.  It is clear that the \pth-power map factors through $\Cal LX^{(p)}$, as $\lambda ^p=\mu ^{(p-1)}\Delta ^{(p-1)}$.  Furthermore, the following diagram of fibrations commutes.
 $$\xymatrix{ \Om X\ar [d]^i\ar [rr]^{\Delta ^{(p-1)}}&&(\Om X)^p\ar [d]^i\ar[rr] ^{\mu^{(p-1)}}&&\Om X\ar [d]^i\\
 \Cal L X\ar [d]^e\ar [rr]^{\Delta ^{(p-1)}}&&\Cal L X^{(p)}\ar [d]^{e^{(p)}}\ar[rr] ^{\mu^{(p-1)}}&&\Cal L X\ar [d]^e\\
X\ar [rr]^{=}&& X\ar[rr] ^{=}&& X}$$

Using techniques similar to those applied in section 2, it is possible to show that there is a twisted bimodule extension $C^*K\tot ( \bc C^*K)^{\otimes p}$ and a quasi-isomorphism 
$$\xymatrix{\Ups ':C^*K\tot (\bc C^*K)^{\otimes p}\ar [r]^(0.6)\simeq &CU^*(\Cal LX^{(p)} ) .}$$ 
This construction can be performed with sufficient naturality to ensure that the diagram
$$\xymatrix{C^*K\tot (\bc C^*K)^{\otimes p} \ar [d]^{\Ups '}\ar [rr]^(0.55){Id\otimes (\psi_K^\sharp)^{( p-1)}}&&C^*K\otimes \bc C^*K\ar [d]^\Ups\\
CU^*(\Cal LX^{(p)})\ar [rr] ^{\Delta ^{(p-1)}}&&CU^*(\Cal LX)}$$
 commutes.  
 
 To complete the construction of a model of the \pth-power map, we need only to find a model of $\mu ^{(p-1)}$.  Modeling $\mu ^{(p-1)}$ is very technical in the general case, however, requiring a fine analysis of the images of $\Ups$ and $\Ups'$.  For certain spaces, we can nevertheless show relatively easily that an acceptable model of $\mu ^{(p-1)}$ is $Id\otimes \chi ^{(p-1)}$, where $\chi$ denotes the usual coproduct on $\bc C^*K$, an esthetically pleasing result.  

We show below that $Id\otimes( \psi _K^\sharp)^{(p-1)}\chi ^{(p-1)}$ is a model of $\lambda ^p$, at least when $K$ is a simplicial model of an odd sphere.  To simplify calculations somewhat, we consider here only the case $p=2$; the case of arbitrary $p$, for a larger class of spaces, can be found in \cite {HR}.   

 We verify first that our candidate to be a model of $\lambda ^2$ is in fact a cochain map.

\proclaim{Proposition 3.1.1}  If $K$ be a finite-type, $1$-reduced simplicial set such that $\psi_K ^{\sharp}$ is commutative, then $\xymatrix@1{Id\otimes (\psi _K^{\sharp}\chi):C^*K\tot\bc C^*K\ar [r]^{}&C^*K\tot \bc C^*K}$ is a cochain map.\endproclaim

\demo{Proof} We need first to show that 
$$(Id\otimes (\psi _K^{\sharp}\chi))\bbd =\bbd(Id\otimes (\psi _K^{\sharp}\chi)).$$
The computation, while combinatorially technical, is not subtle.  The formulas in section 1.4 tell us that if $\psi _K^{\sharp}$ is commutative and $dx_i=0$ for all $i$, then (up to signs)
$$\align
\bbd(1\otimes sx_1|\cdots|sx_n)=&\sum _{i=1}^{n-1}\biggl[\pm \pi(sx_1\star\cdots\star sx_{i-1}\star sx_n)\otimes sx_i|\cdots|sx_{n-1}\\
&\qquad\pm \pi(sx_1\star\cdots\star sx_i)\otimes sx_{i+1}|\cdots|sx_n\\
&\qquad \pm 1\otimes sx_1|\cdots|s(x_ix_{i+1})|\cdots |sx_n\biggr]
\endalign$$
and
$$\align
(a\otimes 1)\cdot (1\otimes sx_1|\cdots|sx_n)=&a\otimes sx_1|\cdots|sx_n\\
&+\sum _{i=1}^{n-1}\biggl[\pm \pi(sa\star sx_1\star\cdots\star sx_{i-1}\star sx_n)\otimes sx_i|\cdots|sx_{n-1}\\
&\qquad\pm \pi(sa\star sx_1\star\cdots\star sx_i)\otimes sx_{i+1}|\cdots|sx_n\biggr],
\endalign$$
while for all $a\in C^*K$
$$\align
(Id\otimes (\psi _K^{\sharp}\chi))(a\otimes sx_1|\cdots|sx_n)=&2a\otimes sx_1|\cdots|sx_n\\
&+a\otimes \sum _{i=1}^{n-1}(sx_1|\cdots|sx_i)\star (sx_{i+1}|\cdots|sx_n).\endalign$$
A bit of elementary algebra and careful counting enable us to show that $Id\otimes (\psi _K^{\sharp}\chi)$ is indeed differential, using the formulas above.
\qed\enddemo

 \proclaim {Theorem 3.1.2}  Let $K$ be the simplicial model of $S^{2n+1}$ with exactly two nondegenerate simplices, in degrees $0$ and $2n+1$, where $n>0$.   The diagram
$$\xymatrix{C^*K\tot\bc C^*K\ar [rr]^{Id\otimes (\psi _K^{\sharp}\chi)}\ar [d]_\simeq^\Ups&&C^*K\tot \bc C^*K\ar [d]_\simeq^\Ups\\
CU^*\Cal LS^{2n+1}\ar [rr]^{CU^*\lambda ^2}&&CU^*\Cal LS^{2n+1}}\tag 3.1.1$$
commutes up to cochain homotopy.
\endproclaim

\demo{Proof}
Recall that $\H^*(\Cal LS^{2n+1})=\Lambda x\otimes \Gamma y$, where $|x|=2n+1$ and $|y|=2n$.  Furthermore $\H ^*(\Om S^{2n+1})=\Gamma y$, and the restriction of $\lambda ^2$ to $\Om S^{2n+1}$ induces an endomorphism of $\Gamma y$  specified by 
$$\align
\H ^*(\lambda ^2|_{\Om S^{2n+1}})(y(m))=&\H^*\chi \H^*\mu  (y(m))\\
=&\H ^*\chi (\sum_{i=0}^m y(i)\otimes y(m-i))\\
=&\sum_{i=0}^m y(i)\star y(m-i)\\
=&2^m y(m).\endalign$$

Consequently, since $\lambda ^2 \circ i=i\circ \lambda ^2|_{\Om S^{2n+1}}$, the endomorphism $\H ^*\lambda ^2$ of $\Lambda x\otimes \Gamma y$ induced by $\lambda ^2$ must satisfy $\H ^*\lambda ^2(1\otimes y(m))=2^m\otimes y(m)$, for degree reasons.  Because $\H ^*\lambda ^2$ is a map of algebras, it is therefore true that
$$\H^*\lambda ^2(x\otimes y(m))= 2^m \cdot x\otimes y(m).$$

Let $z$ denote the unique nondegenerate simplex of $K$ in degree $2n+1$.  As seen in section 2.4, the quasi-isomorphism $$\xymatrix@1{\Ups: C^*K\tot \bc C^*K=(\Lambda z\otimes\Gamma sz,0) \ar [r]^(0.65)\simeq &CU^*(\Cal LS^{2n+1})}$$
sends $z$ to a representative $\zeta$ of $x$, $sz$ to $\xi _1=\omega _0(\zeta)$, which represents $y$, and $sz(m)$ to some representative $\xi _m$  of $y(m)$.  Furthermore, calculations identical to those above show that 
$$(Id\otimes \psi _K^\sharp \chi)(1\otimes sz(m))=2^m \otimes sz(m)\quad \text{and}\quad\text (Id\otimes \psi _K^\sharp \chi)(z\otimes sz(m))=2^m \cdot z\otimes sz(m).$$ Hence, 
$$\Ups (Id\otimes \psi _K^\sharp \chi)(1\otimes sz(m))= 2^m\Ups (1\otimes sz(m)) =2^m \cdot \xi _m,$$
which is a representative of $2^m y(m)$, as is
$$CU^*\lambda ^2\Ups (1\otimes sz(m))=CU^* \lambda ^2 (\xi _m).$$
Since $y(m)$ is the unique class of degree $2nm$, there exists $\varsigma_m\in CU^{2nm-1}\Cal LX$ such that 
$$d^\sharp \varsigma_m =\Ups (Id\otimes \psi _K^\sharp \chi)(1\otimes sz(m))-CU^*\lambda ^2\Ups (1\otimes sz(m)), $$ 
which implies that
$$d^\sharp (\varsigma_m\cdot \zeta) =\Ups (Id\otimes \psi _K^\sharp \chi)(z\otimes sz(m))-CU^*\lambda ^2\Ups (z\otimes sz(m)). $$
Thus, the diagram (3.1.1) commutes up to a cochain homotopy $G$ defined by $G(1\otimes sz(m))=\varsigma _m$ and $G(z\otimes sz(m))=\varsigma _m\cdot \zeta$.
\qed\enddemo

\subhead 3.2 Topological cyclic homology\endsubhead

\noindent As explained in the Preface, we can now construct a cochain complex $tc ^*(X)$ such that $\H^*(tc^*(X)\otimes \Bbb F_2)\cong \H^*((TC(X;2);\Bbb F_2)$, the mod 2 spectrum cohomology of $TC(X;2)$, at least for certain spaces $X$.  The model $tc^*(X)$ is the mapping cone of the following composition, where $\pi$ denotes the obvious projection map. 
$$\xymatrix{(\Lambda \upsilon\otimes C^*K\tot \bc C^*K, \tD)\ar @{->>}[r]^(0.6)\pi\ar [dr]_{\tilde \pi}&C^*K\tot \bc C^*K\ar [d]^{Id\otimes (Id-\psi _K^{\sharp}\chi)}\\
&C^*K\tot\bc C^*K}$$

Recall that the mapping cone of a cochain map $\xymatrix@1{f:(V,d)\ar [r]&(V',d')}$ is a cochain complex $C_f:=(V'\oplus sV, D_f)$, where $D_f(v')=d'v'$ for all $v'\in V'$ and $D_f(sv)=f(v)-s(dv)$ for all $sv\in sV$. It is an easy exercise to show that if $f$ is cochain homotopic to $g$, then $C_f$ and $C_g$ are cochain equivalent. Theorem 3.1.1 suffices therefore to ensure that the mapping cone of the composition above has the right cohomology.

The next theorem now follows immediately from the results of the preceding chapters and section, according the justfication in \cite {HR} of our method of construction of $tc^*(X)$.

\proclaim{Theorem 3.2.1}  Let $X$ be a $1$-connected space with the homotopy type of a finite-type CW-complex, and let $\xymatrix@1{\Ups: C^*K\tot \bc C^*K\ar [r]^(0.6)\simeq & CU^*\Cal LX}$ be a thin free loop model for $X$ such that all primitives of $\bc C^*K$ are indecomposable  and such that
$$\xymatrix{C^*K\tot\bc C^*K\ar [rr]^{Id\otimes (\psi _K^{\sharp}\chi)}\ar [d]_\simeq^\Ups&&C^*K\tot \bc C^*K\ar [d]_\simeq^\Ups\\
CU^*\Cal LX\ar [rr]^{CU^*\lambda ^2}&&CU^*\Cal LX}$$
commutes. Suppose that $\Cal L X_{hS^1}$ has a thin model
$$\xymatrix{{\wU}:(\Lambda \upsilon\otimes C^*K\tot \bc C^*K, \tD)\ar [r]^(0.55)\simeq &(\Lambda \upsilon\otimes CU^*\Cal LX,D^\sharp)}.$$ 
Let
$$tc^*(X)=(C^*K\tot \bc C^*K\oplus s(\Lambda \upsilon\otimes C^*K\tot \bc C^*K), D_{\tilde\pi})$$
where for all $x\otimes c\in C^*K\otimes \bc C^*K$,
$$D_{\tilde \pi }(x\otimes c)= \bbd (x\otimes c)$$
while 
$$\align
D_{\tilde \pi}\bigl(s(1\otimes x\otimes c)\bigr )=&x\otimes c-x\otimes \psi _K^{\sharp}\chi (c)\\
&-s\bigl (1\otimes \bbd(x\otimes c) +\upsilon \otimes S(x\otimes c)\bigr)
\endalign$$
and for $k>0$,
$$D_{\tilde \pi}\bigl(s(\upsilon ^k\otimes x\otimes c)\bigr )=-s\bigl (\upsilon^k\otimes \bbd(x\otimes c) +\upsilon^{k+1} \otimes S(x\otimes c)\bigr).$$
Then $\H^*(tc^*(X)\otimes \Bbb F_2)$ is isomorphic to the mod 2 spectrum cohomology of $TC(X;2)$
\endproclaim

A more general version of this theorem will appear in \cite {HR}.

We conclude this chapter and this article with an example illustrating the use of the models we have built.

\example{Example} Let $n>0$, and let $K$ be the model of $S^{2n+1}$ with exactly one nondegenerate simplex $z$ of positive dimension, in dimension $2n+1$.  As explained in the proof of Theorem 2.4.3,  $C^*K=\Lambda z$, with trivial differential. Furthermore, $\bc C^*K$ is isomorphic as an algebra to $\Gamma sz$, for degree reasons,   and $$fls^*(S^{2n+1})=(\Lambda z\otimes \Gamma sz, 0).$$ 
Thus, $\H^*(\Cal LS^{2n+1})\cong \Lambda z\otimes \Gamma sz$, as has long been known.

Let $sz(m)=sz|\cdots|sz\in \perp^m sz$. We then have
$$hos^*(S^{2n+1})=(\Lambda \upsilon\otimes \Lambda z\otimes \Gamma sz,\tD)$$ 
where $\tD(\upsilon ^k\otimes 1\otimes sz(m))=0$ for all $k$ and $m$, while $$\tD(\upsilon ^k\otimes z\otimes sz(m))=(m+1)\upsilon ^{k+1}\otimes 1\otimes sz(m+1).$$
The integral cohomology of the homotopy orbit space is therefore
$$\H^*(\Cal LS^{2n+1}_{hS^1})\cong \Lambda \upsilon\oplus \Gamma sz\oplus \bigoplus _{k,m\geq 1}\Bbb Z/m\Bbb Z\cdot (v^k\otimes sz(m))$$
as graded modules, while its mod $p$ cohomology is of the form
$$\align
\H^*(\Cal LS^{2n+1}_{hS^1};\Bbb F_p)\cong& \Lambda \upsilon\otimes \Bbb F_p\oplus \Gamma sz\otimes \Bbb F_p\\
&\oplus \bigoplus \Sb{k,m\geq 1}\\ p\mid m\endSb\Bbb F_p\cdot (v^k\otimes sz(m))\oplus \bigoplus\Sb k\geq 0\\ p\mid m+1\endSb \Bbb F_p \cdot (v^k\otimes z\otimes sz(m))\endalign$$
where a Bockstein sends the class of $v^k\otimes z\otimes sz(m)$ to the class of $v^{k+1}\otimes sz(m+1)$.

Finally
$$tc^*(S^{2n+1})=(\Lambda z\otimes \Gamma sz\oplus s(\Lambda \upsilon\otimes \Lambda z\otimes \Gamma sz), D_{\tilde \pi})$$
where for all $m$,
$$D_{\tilde \pi }(z\otimes sz(m))= 0=D_{\tilde \pi }(1\otimes sz(m))$$
while 
$$\align
D_{\tilde \pi}\bigl(s(1\otimes z\otimes sz(m))\bigr )=&z\otimes sz(m)-z\otimes \psi _K^{\sharp}\chi (sz(m))\\
&-s\bigl (\upsilon \otimes S(z\otimes sz(m))\bigr)\\
=&(1-2^m)z\otimes sz(m)-(m+1)s\bigl (\upsilon\otimes 1\otimes sz(m+1))
\endalign$$
and 
$$\align
D_{\tilde \pi}\bigl(s(1\otimes 1\otimes sz(m))\bigr )=&1\otimes sz(m)-1\otimes \psi _K^{\sharp}\chi (sz(m))\\
=&(1-2^m)\otimes sz(m)
\endalign$$
and for $k>0$,
$$D_{\tilde \pi}\bigl(s(\upsilon ^k\otimes z\otimes sz(m))\bigr )=-(m+1)s\bigl (\upsilon^{k+1} \otimes 1\otimes sz(m+1))\bigr)$$ 
and 
$$D_{\tilde \pi}\bigl(s(\upsilon ^k\otimes 1\otimes sz(m))\bigr )=0.$$
Note that we use that 
$$\psi _K^\sharp \chi (sz(m))=\psi _K^\sharp (\sum _{i=0}^m sz(i)\otimes sz(m-i))
=\sum _{i=0}^m \binom mi sz(m)
=2^m sz(m).$$

Modulo 2, these formulas become
$$D_{\tilde \pi }(z\otimes sz(m))= 0=D_{\tilde \pi }(1\otimes sz(m))$$
while 
$$D_{\tilde \pi}\bigl(s(1\otimes z\otimes sz(m))\bigr )=\left \{ \alignedat 2 &z\otimes sz(m)&: &\;m \text{ odd}\\
													&z\otimes sz(m)+s(\upsilon \otimes 1\otimes sz(m+1))&: &\; m\text { even}\endaligned\right .$$
and 
$$D_{\tilde \pi}\bigl(s(1\otimes 1\otimes sz(m))\bigr )=1\otimes sz(m)$$
and for $k>0$,
$$D_{\tilde \pi}\bigl(s(\upsilon ^k\otimes z\otimes sz(m))\bigr )=\left \{ \alignedat 2 &0&: &\;m \text{ odd}\\
													&s(\upsilon^{k+1} \otimes 1\otimes sz(m+1))&: &\;m\text { even}\endaligned\right .$$ 
and 
$$D_{\tilde \pi}\bigl(s(\upsilon ^k\otimes 1\otimes sz(m))\bigr )=0.$$
We can now compute easily that 
$$\align
\H^*(TC(S^{2n+1};2);\Bbb F_2)\cong &\bigoplus \Sb k\geq 1\\ m\geq 0\endSb \Bbb F_2\cdot s(v^k\otimes sz(2m))\oplus \Bbb F_2 \cdot s(v^k\otimes x\otimes sx(2m+1))\\
&\oplus \bigoplus _{m\geq 0}\Bbb F_2\cdot s(v\otimes sz(2m+1))\endalign$$
as graded vector spaces.
\endexample



\Refs\nofrills{References}

\widestnumber\key {HPST}

\ref\key AH \by Adams, J.F.; Hilton, P.J.\paper On the chain algebra of a
loop space\jour Comment. Math. Helv. \vol 30\yr 1956\pages
305--330\endref

\ref \key An \by  Anick, D. \pages 417-453   \paper
Hopf algebras up to homotopy  \yr 1989   
\jour J. Amer.  Math. Soc. \vol 2 \endref

\ref \key B \by Baues, H.-J.\paper The cobar construction as a Hopf 
algebra\jour Invent. Math.\vol 132\yr 1998\pages 467--489\endref

\ref\key Be\by Berrick, A.J.\paper Algebraic K-theory and algebraic topology\inbook Contemporary Developments in Algebraic K-Theory (ed. M Karoubi, A O Kuku, C Pedrini)\bookinfo ICTP Lecture Notes 15, The Abdus Salam ICTP (Trieste, 2004)\pages  97-190\endref

\ref\key Bl\by Blanc, S.\paper Mod\`eles tordus d'espaces de lacets libres et fonctionnels\paperinfo Thesis, EPFL\yr 2004\endref

\ref\key BH\by Blanc, S.; Hess, K.\paper Simplicial and algebraic models for the free loop space\paperinfo In preparation\endref

\ref\key BHM
\by B\"okstedt, M.; Hsiang, W.C.; Madsen, I.
\paper The cyclotomic trace and algebraic $K$-theory of spaces
\jour Invent. Math.
\vol 111
\yr 1993
\pages 465--539
\endref

\ref\key BO1\by B\"okstedt, M.; Ottosen, I.\paper A splitting result for the free loop space of spheres and projective spaces\paperinfo arXiv:math.AT/0411594\endref

\ref\key BO2\by B\"okstedt, M.; Ottosen, I.\paper A spectral sequence for string cohomology \paperinfo arXiv:math.AT/0411571 \endref

\ref\key Br \by Brown, E. H.\paper Twisted tensor products\jour Ann. Math.\vol
69\yr 1959\pages 223--242\endref

\ref\key DH1
\by N. Dupont and K. Hess
\paper Noncommutative algebraic models for fiber squares\vol 314
\pages 449--467
\jour Math. Annalen \yr 1999\endref

\ref\key DH2
\by N. Dupont and K. Hess
\paper How to model the free loop space algebraically
\vol 314
\pages 469--490
\jour Math. Annalen \yr 1999\endref

\ref\key DH3
\by N. Dupont and  K. Hess
\paper Commutative free loop space models at large primes
\yr 2003
\jour Math. Z.\vol 244\pages 1-34\endref

\ref\key DH4
\by N. Dupont and K. Hess 
\paper An algebraic model 
for homotopy fibers\jour Homology, Homotopy and Applications\vol 4\yr 2002\pages 
117--139\endref

\ref\key FHT\by F\'elix, Y.; Halperin, S.; Thomas, J.-C.\book Rational Homotopy Theory\bookinfo Graduate Texts in Mathematics \vol 205\publ Springer\yr 2001\endref

\ref\key GJ\by Goerss, P.G.; Jardine, J.F.\book Simplicial Homotopy Theory\publ Birkh\"auser\bookinfo Progress in Mathematics\vol 174\yr 1999\endref

\ref\key GM \by Gugenheim, V.K.A.M; Munkholm, H.J.\paper On the extended 
functoriality of Tor and Cotor\jour J. Pure Appl. Algebra\vol 4\yr 
1974 \pages 9--29\endref

\ref\key H\by Hess, K.\paper Model categories in algebraic topology\jour Applied Categorical Structures\vol 10\yr 2002\pages 195-220\endref

\ref\key H2\by Hess, K.\paper Algebraic models of homotopy orbit spaces\paperinfo In preparation\endref 



\ref\key HPST\by Hess, K.; Parent, P.-E.; Scott, J.; Tonks, A.\paper A canonical enriched Adams-Hilton  model for simplicial sets\paperinfo Submitted\pages 29 p.\endref

\ref\key HR\by Hess, K.; Rognes, J.\paper Algebraic models for topological cyclic homology and Whitehead spectra of simply-connected spaces \paperinfo In preparation\endref

\ref\key Ho\by Hovey, M.\book Model Categories\bookinfo Mathematical Survey and Monographs\vol 63 \publ American Mathematical Society\yr 1999\endref

\ref\key J\by Jones, J.D.S.\paper Cyclic homology and equivariant homology\jour Invent. Math.\vol 87\yr 1987\pages 403-423\endref

\ref\key K\by Kuribayashi, K.\paper The cohomology of a pull-back on $\Bbb K$-formal spaces\jour Topology Appl.\vol 125\yr 2002\pages 125-159\endref

\ref\key KY\by Kuribayashi, K.; Yamaguchi, T.\paper The cohomology algebra of certain free loop spaces\jour Fund. Math.\vol 154\yr 1997\pages 57-73\endref

\ref\key Man\by Mandell, M. \paper $E_{\infty}$-algebras and $p$-adic 
homotopy theory\jour Topology\vol 40\yr 2001\pages 43--94\endref

\ref\key Mas\by Massey, W.\book Singular Homology Theory\publ Springer\yr 1980\bookinfo Graduate Texts in Mathematics\vol 70\endref 

\ref\key May \by May, J.P.\book Simplicial Objects in Algebraic Topology\publ
University of Chicago Press\yr 1967\bookinfo Midway reprint 1982\endref

\ref\key Mc\by McCleary, J.\book A User's Guide to Spectral Sequences, Second Edition\bookinfo Cambridge studies in advanced mathematics\vol 58 \yr 2001\publ Cambridge University Press\endref

\ref\key Me\by Menichi, L.\paper On the cohomology algebra of a 
fiber\jour Algebr. Geom. Topol. \vol 1 \yr 2001\pages 719--742\endref

\ref\key Me2\by Menichi, L. \paper The cohomology ring of free loop spaces\jour Homology Homotopy Appl.\vol 3\yr 2001\pages 193-224\endref

\ref\key Mi\by Milgram, R.J.\paper Iterated loop spaces\jour Ann. of Math.\vol 84\yr 1966\pages 386--403\endref

\ref\key NT\by Ndombol, B.; Thomas, J.-C.\paper On the cohomology algebra of free loop spaces\jour Topology\vol 41\yr 2002 \pages 85--106\endref

\ref\key R\by Rognes, J.\paper The smooth Whitehead spectrum of a point at odd regular primes\jour Geom. Topol. \vol 7 \yr2003\pages 55--184\endref

\ref\key Sm\by Smith, L.\paper The Eilenberg-Moore spectral sequence and the mod $2$ cohomology of certain free loop spaces\jour Ill. J. Math.\vol 28\yr 1984\pages 516--522\endref

\ref \key SV \by  Sullivan, D.;Vigu\'e-Poirrier, M. \pages 633-644  
\paper The homology theory of the closed geodesic problem \yr  1976  
\jour J. Diff. Geometry \vol 11 \endref

\ref\key S \by Szczarba, R.H.\paper The homology of twisted cartesian 
products\jour Trans. Amer. Math. Soc.\vol 100\yr 1961\pages 
197--216\endref

\ref \key W1 \by Waldhausen, F.
\paper Algebraic K-theory of spaces
\inbook Algebraic and geometric topology
\procinfo Proc. Conf., New Brunswick/USA 1983
\bookinfo Lect. Notes Math. \vol 1126 \pages 318--419 \yr 1985
\endref

\ref \key W2 \by Waldhausen, F.; Jahren, B. ; Rognes, J.
\paper The stable parametrized h-cobordism theorem
\paperinfo In preparation
\endref



\endRefs
\enddocument
\end